\documentclass{article}%
\usepackage{amssymb}
\usepackage{amsfonts}
\usepackage{amsmath}
\usepackage{graphicx}%
\providecommand{\U}[1]{\protect\rule{.1in}{.1in}}
\newtheorem{theorem}{Theorem}

\newtheorem{corollary}[theorem]{Corollary}

\newtheorem{lemma}[theorem]{Lemma}

\newtheorem{proposition}[theorem]{Proposition}

\begin{document}

\title{Hartogs extension for systems of differential equations}
\author{V. P. Palamodov\\Tel Aviv University}
\date{}
\maketitle

\textbf{Abstract.} The phenomenon of removable singularity is studied for
overdetermined elliptic systems of differential equations. We show that the
dimension of the characteristic variety of the system plays a key role in the problem.

\textbf{MSC} 2010\ \ 58J10, 35N10, 35A27

\section{Introduction}

According to the famous theorems of Hartogs and Osgood-Brown any compact
singularity (with no holes) of a holomorphic function of several variables is
removable. This fact can be viewed as a property of solutions of the
Cauchy-Riemann system of differential equations\ in a domain of the real space
$\mathbb{R}^{2m}$ where $m>1.$ This phenomenon was extended by Ehrenpreis
\cite{Eh} who stated that a compact singularity is always removable for any
system of two equations with constant coefficients and relatively prime
symbols. It was shown in \cite{P1} that the automatic extension property of
solutions of a general system $\mathit{M}$ of equations with constant
coefficients is governed by the modules $\mathrm{Ext}^{k}\left(
\mathit{M,D}\right)  ,$ $k=1,2,...$(see below). In particular the equation
$\mathrm{Ext}^{1}=0$ guarantees absence of compact singularity. Vanishing of
higher $\mathrm{Ext}$ implies removing of some noncompact singularities in
particular for singularities supported by submanifolds. Kawai \cite{K} and
Kawai and Takei \cite{KT} treated systems of equations with one unknown
function and commuting operators. They stated automatic extension of solutions
in the class of hyperfunctions. According to \cite{K} (Theorem 2) the last
property is however weaker than automatic extensibility in the classical sense.

We study here the phenomenon of compulsory extension for systems of linear
partial differential equations with analytic coefficients in $\mathbb{R}^{n}$
of the general form%
\begin{equation}
P\left(  x,\partial_{x}\right)  u=0 \label{1}%
\end{equation}
Here $\partial_{x}\doteq\left(  \partial/\partial x_{1},...,\partial/\partial
x_{n}\right)  ,\ P=\{p_{ij}\}$ is a $s\times r$-matrix differential operator,
$u=\left(  u_{1},...,u_{r}\right)  $ are unknown functions, numbers $s$ and
$r$ are arbitrary. We show that for any elliptic system a small compact
singularity is always removable (in the classical sense) if dimension $d$ of
the characteristic variety $\mathit{V}$ is strictly less than $n-1$ just as
for the case of Cauchy-Riemann system with $m>1$. In the case of smaller $d$ a
more strong statement holds (Theorem \ref{ee}). In particular, an analytic
submanifold $S$ of dimension $s$ can not be a support of a nonremovable
singularity of a solution if $s<n-d-1$. This is not the case if $s=n-d-1.$

\section{Regularity conditions for a differential matrix}

We do not assume any special structure of the matrix $P$ but impose a general
condition of regularity which is well known in several forms \cite{R}%
,\cite{M2}. Fix an arbitrary point $x\in X$, denote by $\mathit{O}_{x}$ the
algebra of germs of analytic functions at the point $x\in\mathbb{R}^{n}$ and
by $\mathit{D}_{x}$ the algebra of differential operators in $\mathit{O}_{x}$.
Let $P_{i}=\left(  p_{i1},...,p_{ir}\right)  ,$ $i=1,...,s$ be rows of the
matrix $P;$ consider a linear combination
\begin{equation}
Q\left(  x,\partial_{x}\right)  =\sum_{i=1}^{s}a_{i}\left(  x,\partial
_{x}\right)  P_{i}\left(  x,\partial_{x}\right)  \in\mathit{D}_{x}%
^{r}\label{2}%
\end{equation}
with some $a_{i}\in\mathit{D}_{x}$ where $Q=\left(  q_{1},...,q_{r}\right)
$We assume that there exist $b_{i}\in\mathit{D}_{x},i=1,...,s$ such that%
\begin{equation}
Q\left(  x,\partial_{x}\right)  =\sum b_{i}\left(  x,\partial_{x}\right)
P_{i}\left(  x,\partial_{x}\right)  \label{3}%
\end{equation}
and $\mathrm{ord\,}b_{i}+\mathrm{ord\,\,}p_{ij}\leq\mathrm{ord\,}q_{j}$ for
all $i=1,...,s$, $j=1,...,r$ (\textrm{ord}$a$ means the order of a
differential operator $a$). In other words, there is no cancellations of
higher order terms in the right-hand side of (\ref{3}). This condition is not
in fact restrictive since it can be always satisfied if the matrix $P$ is
supplemented by several lines of the form (\ref{2}).\newline%
\textbf{Definition. }Fix some integers $\sigma_{1},...,\sigma_{s}$ and
$\rho_{1},...,\rho_{t}$ (called shifts) such that
\[
\mathrm{\deg\ }p_{ij}\leq\sigma_{i}-\rho_{j},\ i=1,...,s;~j=1,...,r
\]
The principal part of the system (\`{a} la Douglis-Nirenberg \cite{DN}) is the
matrix $\mathrm{P}=\{\mathrm{p}_{ij}\}$, where $\mathrm{p}_{ij}$ is the sum of
homogeneous terms of $p_{ij}$ of degree $\sigma_{i}-\rho_{j}$ ($p_{ij}=0$ if
there is no such terms). Substituting partial derivatives $\partial/\partial
x_{i}$ by independent variables $\xi_{i},i=1,...,n,$ we\ obtain homogeneous
polynomials $\mathrm{p}_{ij}\left(  x,\xi\right)  \mathit{\ }$in $\xi=\left(
\xi_{1},...,\xi_{n}\right)  \ $with coefficients in $\mathit{O}_{x}.$ The next
condition is essential: \newline(*) For any point $x\in X$ and
polynomials$\ \mathrm{r}_{1},...,\mathrm{r}_{s}\in\mathbb{C}[\xi_{1}%
,...,\xi_{n}]$ such that%
\[
\sum_{i}\mathrm{r}_{i}\left(  \xi\right)  \mathrm{P}_{i}\left(  x,\xi\right)
=0
\]
where \textrm{P}$_{i}\left(  x,\xi\right)  $ denotes the vector$\ \left(
\mathrm{p}_{i1}\left(  x,\xi\right)  ,...,\mathrm{p}_{ir}\left(  x,\xi\right)
\right)  $ there exist functions$\ $\newline$\mathrm{R}_{1},...,\mathrm{R}%
_{s}$ $\in\mathit{O}_{x}[\xi_{1},...,\xi_{n}]$ such that%
\[
\sum_{i}\mathrm{R}_{i}\left(  y,\xi\right)  \mathrm{P}_{i}\left(
y,\xi\right)  =0
\]
for $y\ $in a neighborhood of the point $x$ such that and $\mathrm{R}%
_{i}\left(  x,\xi\right)  =\mathrm{r}_{i}\left(  \xi\right)  ,i=1,...,s.$ In
fact, this condition need to be checked for only finite number of vectors
$\left(  \mathrm{r}_{1},...,\mathrm{r}_{s}\right)  $ and it is generic that is
(**) is always fulfilled in the compliment to a nowhere dense analytic set
\cite{P2}.

Note that in the case $r=s=1$ the condition (*) only means that the principal
part \textrm{P\ }of $P$ does not vanish at $x.$

\section{Differential modules and filtrations}

Now we rewrite the above conditions in a more algebraic form. Let again
$x\in\mathbb{R}^{n}\ $and $\mathit{D}$ be the algebra of differential
operators in $\mathbb{R}^{n}$ with coefficients in the algebra $\mathit{O}$ of
germs at $x$ of analytic functions in $\mathbb{R}^{n}\ $(here and later we
omit the subscript $x$). The algebra $\mathit{D}$ has natural filtration
$\{\mathit{D}_{k},\ k=0,1,...\},$ where $\mathit{D}_{k}$ is the $\mathit{O}%
$-module of differential operators $a\in\mathit{D}$ of order
\textrm{ord\thinspace}$a\leq k$ and $\mathit{D}_{0}=\mathit{O.}$ The
associated graded module
\[
\mathrm{D}=\mathrm{gr\,}\mathit{D}=\oplus_{k=0}^{\infty}\mathit{D}%
_{k}/\mathit{D}_{k-1}%
\]
is a commutative $\mathit{O}$-algebra. Fix a coordinate system $x_{1}%
,...,x_{n}$ in $\mathbb{R}^{n}.$ The algebra $\mathrm{D}$ is isomorphic to the
graded algebras $\mathit{O}[\xi_{1},...,\xi_{n}]$ where the generator $\xi
_{i}$ is represented by the operator $\partial/\partial x_{i},i=1,...,n.$ The
algebra $\mathrm{D}\otimes_{\mathit{O}}\mathbb{C}$ is then isomorphic to the
graded algebra of homogeneous polynomials in $T_{x}^{\ast}\left(
\mathbb{R}^{n}\right)  .$

Fix a natural $r$ and a vector $\rho=\left(  \rho_{1},...,\rho_{r}\right)
\in\mathbb{Z}^{r}$; the increasing sequence of $\mathit{O}$-submodules
\[
\mathit{D}_{k}^{\rho}=\{a\in\mathit{D}^{r},\mathrm{ord}_{\rho}a\leq
k\},\ k\in\mathbb{Z}%
\]
is called filtration generated by the shift vector $\rho,$ where
$\mathrm{ord}_{\rho}a=\mathrm{ord}_{\rho}\left(  a_{1},...,a_{r}\right)
=\max_{i}\mathrm{ord}a_{i}+\rho_{i}.$ The graded vector space%
\[
\mathrm{D}^{\rho}=\oplus_{k}\mathit{D}_{k}^{\rho}/\mathit{D}_{k-1}^{\rho}%
\]
is a module over the graded commutative algebra $\mathrm{D}$. Let $r,s$ be
natural numbers; any morphism of left $\mathit{D}$-modules $P:\mathit{D}%
^{s}\rightarrow\mathit{D}^{r}$ can be written in the form $a\mapsto aP,$ where
an element $a=\left(  a_{1},...,a_{s}\right)  \in\mathit{D}^{s}$ is thought as
row and $P$ as a $s\times r$-matrix whose entries $p_{ij},i=1,...,s,j=1,...,r$
are sections of $\mathit{D.}$ Let $\sigma$ denote the filtration in
$\mathit{D}^{s}$ defined by a shift vector $\sigma=\left(  \sigma
_{1},...,\sigma_{s}\right)  .$ The morphism $P$ agrees with the filtrations,
if $\mathrm{ord}_{\rho}\left(  aP\right)  \leq\mathrm{ord}_{\sigma}a$ for any
$a\in\mathit{D}^{s}.$ This condition is equivalent to the inequalities
$\mathrm{ord\,}p_{ij}\leq\sigma_{i}-\rho_{j}.$ Let $\mathrm{p}_{ij}$ be the
some of terms of $p_{ij}$ of order $\sigma_{i}-\rho_{j}.$ The matrix
$\mathrm{P}=\{\mathrm{p}_{ij}\}$ is called the principal part of $P$ with
respect to the filtrations generated by $\rho$ and $\sigma.$ The operator $P$
is called \textit{elliptic} at a point $x$ if \textrm{rank }$\mathrm{P}\left(
x,\xi\right)  =s$ for any real $\xi\neq0.$

Let $\mathit{M}$ be a left\textit{\ }$\mathit{D}$\textit{-}module; suppose
that $\mathit{M}$ has a increasing filtration by $\mathit{O}$-submodules
$\mathit{M}_{k},k\in\mathbb{Z\ }$such that$\ \cup\mathit{M}_{k}=\mathit{M}$
and $\mathit{D}_{i}\,\mathit{M}_{k}\subset\mathit{M}_{k+i}$ for any $i$ and
$k.$ Then we call $\mathit{M}$ filtered $\mathit{D}$-module.\textit{ }For such
a module the direct sum
\[
\mathrm{gr\,}\mathit{M}=\oplus_{k=-\infty}^{\infty}\mathit{M}_{k}%
/\mathit{M}_{k-1}%
\]
has a natural structure of $\mathrm{D}$-module.

\textbf{Definition. }Let $\mathit{M}$ and $\mathit{N}$ be filtered left (or
right) $\mathit{D}$-modules. We say that a $\mathit{D}$-morphism $\alpha:$
$\mathit{M}\rightarrow\mathit{N}$ agrees with the filtrations, if
$\alpha\left(  \mathit{M}_{k}\right)  \subset\mathit{N}_{k}$ for
$k=0,1,2,...$. If $\alpha$ agrees with filtrations it generates a morphism of
graded modules $\mathrm{gr\,}\alpha:\mathrm{gr\,}\mathit{M}\rightarrow
\mathrm{gr\,}\mathit{N}$ and the correspondence $\alpha\mapsto\mathrm{gr\,}%
\alpha$ is a functor.

\section{Complexes and symbols}

Let
\begin{equation}
...\rightarrow\mathit{D}^{t}\overset{Q}{\rightarrow}\mathit{D}^{s}\overset
{P}{\rightarrow}\mathit{D}^{r} \label{35}%
\end{equation}
be an exact sequence of left $\mathit{D}$-modules. The morphism $P$ in
(\ref{21}) acts by right multiplication $a\mapsto aP\ $where $P\ $is a
$s\times r$-matrix\ whose entrees $p_{ij}\ $belong to $\mathit{D.\ }$The
morphisms $...,Q$ can be realized in a similar way.\ 

\textbf{Definition. }We say that the modules in (\ref{35}) are supplied with
filtrations generated by some shift vectors $\rho,\sigma,\tau,...$. We shall
denote these modules by $...,\mathit{D}^{\tau},\mathit{D}^{\sigma}%
,\mathit{D}^{\rho}\ $to fix the respective filtrations generated by some shift
vectors $...,\tau,\sigma,\rho$ where in particular $\rho=\left(  \rho
_{1},...,\rho_{r}\right)  ,\ $\textrm{ord}$_{\rho}a=\max_{j}\mathrm{ord~}%
a_{j}+\rho_{j}$.\ We assume that the morphisms $...,Q,P$ agree with these
filtrations which means $...,\mathrm{ord\ }q_{jk}\leq\tau_{j}-\sigma
_{k},\mathrm{ord\,}p_{ij}\leq\sigma_{i}-\rho_{j}$ for entrees of these
matrices. The sequence of graded $\mathrm{D}$-modules is then well-defined%
\[
...\rightarrow\mathrm{D}^{\tau}\overset{\mathrm{Q}}{\rightarrow}%
\mathrm{D}^{\sigma}\overset{\mathrm{P}}{\rightarrow}\mathrm{D}^{\rho}%
\]
By tensoring over the algebra $\mathit{O}$ we get the complex%
\[
...\rightarrow\mathrm{D}^{\tau}\otimes\mathbb{C}\overset{\mathrm{Q}%
\otimes\mathbb{C}}{\rightarrow}\mathrm{D}^{\sigma}\otimes\mathbb{C}%
\overset{\mathrm{P}\otimes\mathbb{C}}{\rightarrow}\mathrm{D}^{\rho}%
\otimes\mathbb{C}%
\]
of free graded modules over the commutative algebra $\mathrm{A}=\mathrm{D}%
\otimes\mathbb{C\cong C}\left[  \xi_{1},...,\xi_{n}\right]  $.$\ $The set of
maximal ideals in the algebra \textrm{A }is isomorphic to $\mathbb{C}^{n}.$
For a maximal ideal $\mathfrak{m}$ in \textrm{A} we take tensor product with
the quotient algebra. This yields a\ complex of $\mathbb{C}$-vector spaces%
\[
...\rightarrow\left(  \mathrm{D}^{\tau}\otimes\mathbb{C}\right)
\mathbb{\otimes}_{\mathrm{A}}\mathrm{A}/\mathfrak{m}\overset{\mathrm{Q}%
_{\mathfrak{m}}}{\rightarrow}\left(  \mathrm{D}^{\sigma}\otimes\mathbb{C}%
\right)  \mathbb{\otimes}_{\mathrm{A}}\mathrm{A}/\mathfrak{m}\overset
{\mathrm{P}_{\mathfrak{m}}}{\rightarrow}\left(  \mathrm{D}^{\rho}%
\otimes\mathbb{C}\right)  \mathbb{\otimes}_{\mathrm{A}}\mathrm{A}/\mathfrak{m}%
\]
\textrm{ }where\textrm{ }$...,\mathrm{P}_{\mathfrak{m}}=\mathrm{P}%
\otimes\mathbb{C\otimes}\mathrm{A}/\mathfrak{m.}$ Because of $\mathrm{A}%
/\mathfrak{m}\cong\mathbb{C}$ this complex can be written in a simple form%
\begin{equation}
...\rightarrow\mathbb{C}^{t}\overset{\mathrm{Q}\left(  x,\xi\right)
}{\rightarrow}\mathbb{C}^{s}\overset{\mathrm{P}\left(  x,\xi\right)
}{\rightarrow}\mathbb{C}^{r}\label{6}%
\end{equation}
where $\xi$ is the point in $\mathbb{C}^{n}$ corresponding to the ideal
$\mathfrak{m.}$ Here $...,\mathrm{Q}\left(  x,\xi\right)  ,\mathrm{P}\left(
x,\xi\right)  $ are matrices whose entrees are analytic functions of $x$ and
polynomial function of $\xi.$

\textbf{Definition. }We call (\ref{6}) the \textit{principal symbol }of
(\ref{35}). The complex (\ref{35}) is called elliptic if the symbol is exact
at any real point $\xi\neq0.$

\section{Local solvability}

Let$\ x\in\mathbb{R}^{n},\ \mathit{D}=\mathit{D}_{x}$ and
\begin{equation}
...\rightarrow\mathit{D}^{t}\overset{Q}{\rightarrow}\mathit{D}^{s}\overset
{P}{\rightarrow}\mathit{D}^{r} \label{4}%
\end{equation}
be an exact sequence of left $\mathit{D}$-modules. Denote by $\mathit{E}$ the
space of germs at $x$ of $C^{\infty}$-functions defined in $\mathbb{R}^{n}.$
This space has the natural structure of a left $\mathit{D}$-module. Applying
the functor $\mathrm{Hom}_{\mathit{D}}\left(  \cdot,\mathit{E}\right)  $ to
(\ref{4}) we obtain a complex of vector spaces:%
\begin{equation}
\mathit{E}^{r}\overset{P}{\rightarrow}\mathit{E}^{s}\overset{Q}{\rightarrow
}\mathit{E}^{t}\rightarrow... \label{5}%
\end{equation}
where the matrices $P,Q,...$ act by left multiplication like in (\ref{1}).

\begin{theorem}
If (\ref{4}) is exact and elliptic, then the sequence (\ref{5}) is exact.
\end{theorem}

The case of arbitrary operator $P$ with constant coefficients is considered in
\cite{M1},\cite{P1}. For the case of analytic coefficients this statement is
essentially due to Malgrange \cite{Ma2} who proved it for
Newlander-Nirenberg's operator by reduction to the case of germs of analytic
functions (the method coming from the Hodge theory). A proof in the general
case was done by Andreotti-Nacinovich \cite{AN} by the same method. Quite
different method was used by H\"{o}rmander \cite{H2}.

We state here a quantitative version of this Theorem.\ Let $\mathit{E}$ be the
sheaf of germs of $C^{\infty}$-functions in $\mathbb{R}^{n}$. The space
$\mathit{E}\left(  U\right)  =\Gamma\left(  U,\mathit{E}\right)  $ for any
open $U\subset\mathbb{R}^{n}\ $has the natural Fr\'{e}chet topology. Any
differential $s\times r$-matrix $P$ as in (\ref{1}) defines for any open
$U\subset X$ a linear continuous operator $P:\mathit{E}^{r}\left(  U\right)
\rightarrow\mathit{E}^{s}\left(  U\right)  $. We denote by $\mathit{E}%
_{P}\left(  U\right)  $ its kernel.

Fix an Euclidean structure in $\mathbb{R}^{n};$ for a point $x\in
\mathbb{R}^{n}$ and a number $r>0$ the notation $U_{x}\left(  r\right)  $
means $r$-neighborhood of $x.$

\begin{theorem}
\label{cs}If (\ref{4}) is elliptic, then\ \newline\textbf{A}. there exists a
continuous function $b_{x}$ in $X$ such that for arbitrary point $x\in X,$
arbitrary $0<r\leq1$ and arbitrary $g\in\mathit{E}_{Q}\left(  U_{x}\left(
r\right)  \right)  $ there exists a $f\in\mathit{E}^{r}\left(  U_{x}\left(
b_{x}r\right)  \right)  $ such that
\begin{equation}
Pf=g \label{25}%
\end{equation}
in $U_{x}\left(  b_{x}r\right)  .$\newline\textbf{B}. There exists a linear
continuous operator $\mathrm{s}_{x,r}:\mathit{E}_{Q}\left(  U_{x}\left(
r\right)  \right)  \rightarrow\mathit{E}^{r}\left(  U_{x}\left(
b_{x}r\right)  \right)  $ that provides a solution to (\ref{25}).\newline
\end{theorem}

$\blacktriangleleft\ $\textbf{1. }We will construct a Laplace-like operator
$\Omega$ for (\ref{4}) and reduce the statement to the case when $\Omega g=0.$
Because of the morphisms agree with the filtrations the inequalities are
fulfilled
\[
\mathrm{ord\,}q_{ij}\leq\tau_{i}-\sigma_{j},\ \mathrm{ord\ }p_{ij}\leq
\sigma_{i}-\rho_{j}%
\]
where $Q=\{q_{ij}\},P=\{p_{ij}\}$ are the entrees of matrices $Q$ and
$P$.\ For any differential operator $a$ in $X\subset\mathbb{R}^{n}$ with
analytic coefficients the formal adjoint operator $a^{\ast}$ acts on smooth
densities and on distributions with compact support in $X$:%
\[
\int_{X}a^{\ast}\left(  v\right)  u=\int_{X}va\left(  u\right)
\]
Identifying a function $u$ with the density $u\mathrm{d}x$ where $\mathrm{d}x$
is the Euclidean volume form in $\mathbb{R}^{n}$ we make the adjoint operator
$a^{\ast}$ acting on functions. It also has analytic coefficients. The
transformation $a\mapsto a^{\ast}$ is $\mathbb{C}$-linear and we have $\left(
ab\right)  ^{\ast}=b^{\ast}a^{\ast}$. The Laplace operator $\Delta$ in
$\mathbb{R}^{n}$ is self-adjoint; we denote$\ \square=-\Delta$. For a natural
$k$ and a vector $\omega=\left(  \omega_{1},...,\omega_{k}\right)  $ with
natural coordinates we denote by $\square^{\omega}$ the diagonal $k\times
k$-matrix $\left(  \square^{\omega_{1}},...,\square^{\omega_{k}}\right)  $.
Set $t=\mathrm{\max}\left(  \tau_{1},...,\tau_{t}\right)  $ and denote
$t+\rho=\left(  t+\rho_{1},...,t+\rho_{r}\right)  ,$ $t-\tau=...$ The
differential operator
\[
\Omega=P\square^{t+\rho}P^{\ast}+\square^{\sigma}Q^{\ast}\square^{t-\tau
}Q\square^{\sigma}%
\]
is well defined in the sheaf $\mathit{E}^{s}$.

\begin{lemma}
\label{DN}$\Omega$ is an elliptic operator in the sense of Douglis-Nirenberg
in $X$ with the shift vectors equal to $2t+2\sigma=\left(  2t+2\sigma
_{1},...,2t+2\sigma_{s}\right)  .$
\end{lemma}

We postpone a proof of Lemma. Because of $\Omega$ is elliptic, there exists a
countable family of local fundamental solutions $\Phi$ defined in open sets
$U_{\Phi}$ such that $X=\cup U_{\Phi}.$ Take a fundamental solution $\Phi$
(Lemma \ref{EF}), an arbitrary point $x\in U_{\Phi}$ and a number $r>0$ such
that $U_{x}\left(  r\right)  \subset U_{\Phi}.$ Choose a smooth cut function
$e$ with$\mathcal{\ }$support in $U_{x}\left(  r\right)  $ that is equal to 1
in $U_{x}\left(  r/2\right)  .$ Suppose that a function $g\in\mathit{E}%
^{s}\left(  U_{x}\left(  r\right)  \right)  $ fulfils $Qg=0$ and set%
\[
g_{e}=eg,\ h=Q^{\ast}\square^{t-\tau}Q\square^{\sigma}\Phi g_{e},\text{
}f=\square^{t+\rho}P^{\ast}\Phi g_{e}%
\]
We have%
\[
Pf=P\square^{t+\rho}P^{\ast}\Phi g_{e}=\Omega\Phi g_{e}-\square^{\sigma
}Q^{\ast}\square^{t-\tau}Q\square^{\sigma}\Phi g_{e}=g_{e}-\square^{\sigma}h,
\]
Thus the function $f$ is a solution of (\ref{25}) modulo a function
$\square^{\sigma}h.$ On the other hand
\[
\Omega h=\square^{\sigma}Q^{\ast}\square^{t-\tau}Q\square^{\sigma}Q^{\ast
}\square^{t-\tau}Q\square^{\sigma}\Phi g_{e}=\square^{\sigma}Q^{\ast}%
\square^{t-\tau}Q\Omega\Phi g_{e}=\square^{\sigma}Q^{\ast}\square^{t-\tau
}Qg=0
\]
in $U_{x}\left(  r/2\right)  $ since $P^{\ast}Q^{\ast}=0$ and $QP=0.$ It
follows that the function $h$ is analytic in $U_{x}\left(  r/2\right)  $ since
of Petrowsky's result \cite{P}.

\textbf{2. }We show that $h$ has holomorphic continuation in a quantified
neighborhood $Z_{x}$ of the point $x.$

\begin{lemma}
\label{EF}Let $A$ be a $s\times s$-matrix differential operator with analytic
coefficients in a ball $U\subset\mathbb{R}^{n}$ that is elliptic in the sense
of Douglis-Nirenberg. Then there exists a fundamental solution $E=E\left(
x,u\right)  $ defined in $U\times U$ that admits a holomorphic extension
$E\left(  z,w\right)  $ in the domain%
\begin{equation}
Z=\{z=x+\imath y,w=u+\imath v,\ c\left\vert y-v\right\vert \leq\left\vert
x-u\right\vert <r\}\ \label{9}%
\end{equation}
for some positive $c$ and $r.$ Moreover $\phi E$ defines a bounded operator in
$L_{2}\left(  U\right)  \rightarrow L_{2}\left(  \mathbb{R}^{n}\right)  $ for
any test function $\phi.$with support in $U.$
\end{lemma}

$\blacktriangleleft$ According to Douglis-Nirenberg's condition the operator
$A$ defines a map $A:\mathit{D}^{\sigma}\rightarrow\mathit{D}^{\rho}$ of order
$0$ where $\mathit{D}^{\sigma}$ and $\mathit{D}^{\rho}$ are free $\mathit{D}%
$-module of rank $s$ with filtrations defined by some shift vectors $\rho$ and
$\sigma.$\ Let $\mathrm{A}\left(  x,D\right)  $ be the principal part of this
operator that is $\mathrm{A}=\{\mathrm{a}_{ij}\}$ where $\mathrm{a}_{ij}$ is
the homogeneous part of $a_{ij}$ of order $\sigma_{i}-\rho_{j}.$The principal
symbol $\mathrm{A}\left(  x,\xi\right)  $ with respect to the vectors
$\rho,\sigma$ is elliptic for any $x\in X$ that is the scalar operator
$\det\mathrm{A}\left(  x,\xi\right)  $ is elliptic. Set $A_{q}\left(
\partial_{x}\right)  =A\left(  q,\partial_{x}\right)  $ and $B=A_{q}-A$ where
$q$ denotes the center of $U.$ The elliptic operator $A_{q}$ with constant
coefficients possesses a fundamental solution $E_{q}\left(  x,y\right)
=E_{q}\left(  x-y\right)  $ in $\mathbb{R}^{n}$ that has a holomorphic
extension $E_{0}\left(  z\right)  $ to the neighborhood of $\mathbb{R}%
^{n}\backslash\{0\}$ of the form $\{z\in\mathbb{C}^{n};\left\vert
\operatorname{Im}z\right\vert <c_{A}\left\vert \operatorname{Re}z\right\vert
\}$ where the constant $c_{A}$ is determined from the condition $\mathrm{\det
\,A}_{q}\left(  \xi+\imath\eta\right)  \neq0$ for $\left\vert \eta\right\vert
<c_{A}\left\vert \xi\right\vert $ where $A_{q}\left(  \zeta\right)  $ is the
symbol of $A_{q}.$ Such a fundamental solution can be written as the
Fourier-Laplace integral of $A_{q}^{-1}\left(  \zeta\right)  $ taken over a
$n$-cycle in $\mathbb{C}^{n}$ that coincides with $\mathbb{R}^{n}$ up to a
compact subset. We construct a fundamental solution for $A$ by the method of
E. Levi. Choose a number $r$ that is smaller than the radius of $U$ and take a
test function$\ \phi$ in $U$ that is equal to $1$ for $\left\vert
x-q\right\vert \leq r$. Consider the series of operators in $L_{2}\left(
\mathbb{R}^{n}\right)  $%
\begin{equation}
E=E_{q}\sum_{k=0}^{\infty}F^{k},\ F=\phi BE_{q} \label{13}%
\end{equation}
where we set $\phi B=0$ in $\mathbb{R}^{n}\backslash U.$ Let $e_{ij}$ be the
$s\times s$-matrix whose entry equals $1$ on $ij$-place and $0$ otherwise and
$b$ be a differential operator in $\mathbb{R}^{n}$ with constant coefficients
such that \textrm{ord}$\,b\leq\sigma_{i}-\rho_{j}.$ The operator $be_{ij}$
defines a map $\mathit{D}^{\sigma}\rightarrow\mathit{D}^{\rho}$ of order $0$
which implies that the composition $be_{ij}E_{q}$ is a bounded operator in
$L_{2}\left(  \mathbb{R}^{n}\right)  ^{s}.$ We have
\[
B=\sum_{i,j=1}^{r}\beta_{ij}\left(  x\right)  b_{ij}e_{ij}%
\]
where for all $i,j,\ $\textrm{ord\thinspace}$b_{ij}\leq\sigma_{i}-\rho_{j}$
and $\beta_{ij}\left(  z\right)  $ are analytic functions that vanish for
$x=q$, since $B\left(  q,\partial_{x}\right)  =0.$ The norm $N=\left\Vert \phi
BE_{q}\right\Vert $ can be made smaller $1/2$ if we take $r=r_{q}$
sufficiently small and the series (\ref{13}) converges as an operator in
$L_{2}\left(  U\right)  ^{s}.$ Moreover we have%
\[
\left\Vert \psi E\right\Vert <2\left\Vert \psi E_{q}\right\Vert <\infty
\]
for any text function $\psi$ since the kernel $E_{q}$ has weak singularity. It
is easy to check that $AE=\mathrm{id}$ in $U.$ The kernel $E\left(
x,u\right)  $ is real analytic out of the diagonal since the operator $A$ is
elliptic. Moreover it has a holomorphic extension to the domain (\ref{9}) with
$c=c_{A_{q}}/2$ due to H\"{o}rmander \cite{H}, Theorem 5.3.3.
$\blacktriangleright$

\begin{lemma}
For an arbitrary linear differential operator $A$ with analytic coefficients
in $X$ that is elliptic in the sense of Douglis-Nirenberg there exists a
positive continuous function $r_{q}$ in $X$ such that for any$\ q\in X,$ any
$r\leq r_{q}$ and arbitrary solution $f$ of the equation $Af=0$ in the ball
$U_{q}\left(  r\right)  $ has a unique holomorphic extension to the ball
$Z_{q}\left(  c_{q}r\right)  $ in $\mathbb{C}^{n}$ with the same center and
radius $c_{q}r$ where $c_{q}=c_{A_{q}}/2$ is a continuous function in $X.$
\end{lemma}

$\blacktriangleleft$ To prove the statement of Lemma we choose a number
$r^{\prime}<r$ and a cut function $e$ supported in $U_{q}\left(  r\right)  $
that is equal to 1 in $U_{q}\left(  r^{\prime}\right)  $ and evaluate the
function $f$ in $U_{q}\left(  r^{\prime}\right)  $ by the integral%
\[
f\left(  x\right)  =\int_{\mathbb{R}^{n}}A^{\ast}\left(  u,\partial
_{u}\right)  e\left(  u\right)  f\left(  u\right)  E\left(  x,u\right)
\mathrm{d}u
\]
We can move now the point $x$ to an arbitrary point $z=x+\imath y$ such that
$\left\vert y\right\vert <c_{q}r^{\prime}/2.$ The function $E\left(
z,u\right)  $ is holomorphic since the support of the integrand is contained
in $U\left(  r\right)  \backslash U\left(  r^{\prime}\right)  $ hence
$\left\vert x-u\right\vert >r^{\prime}.$ This gives a holomorphic extension of
$f$ in the ball $Z_{q}\left(  c_{q}r/2\right)  .$ $\blacktriangleright$

Thus the function $h$ as above has holomorphic extension to the ball
$Z_{x}\left(  b_{x}r\right)  $ for some continuous function $b_{x}.$ This
construction has the property II with any $r<\mathrm{dist}\left(  x,\partial
U_{\Phi}\right)  $ and a positive continuous function $b_{x}=b_{x}\left(
U_{\Phi}\right)  \leq1$ defined in the domain $U_{\Phi}.$ We set $b_{x}\left(
U_{\Phi}\right)  =0$ Take the maximum
\[
b_{x}\left(  X\right)  =\max_{\Phi}\left\{  \frac{\delta\left(  x\right)
}{\delta\left(  x\right)  +1}b_{x}\left(  U_{\Phi}\right)  ,x\in U_{\Phi
}\right\}
\]
where $\delta\left(  x\right)  \doteq\mathrm{dist}\left(  x,\partial U_{\Phi
}\right)  $ for $x\in U_{\Phi}$ and $\delta\left(  x\right)  =0$ for $x\in
X\backslash U_{\Phi},\ $over the family of all fundamental solutions $\Phi.$
It is well defined for any $x\in X$ and is a continuous function$.$ We have
$b_{x}\left(  X\right)  \leq b_{x}\left(  U_{\Phi}\right)  $ for any $\Phi,$
hence the function $b_{x}\left(  X\right)  $ fulfils\textbf{ A} and \textbf{B}.

\textbf{3. }Now it is sufficient to prove the statements of Theorem \ref{cs}
for the sheaf $\mathit{H}$ of germs of analytic functions in $\mathbb{R}^{n}.$
A construction of a solution to (\ref{25}) in the space of germs of analytic
functions can be done by the method of \cite{P2} that guarantees the
properties \textbf{A}, \textbf{B} in terms of balls $Z_{x}\left(  r\right)  $
in $\mathbb{C}^{n}.$ $\blacktriangleright$

\textbf{4.}$\ $Proof of Lemma \ref{DN}. We have%

\begin{align*}
\mathrm{ord}\left(  P\square^{t+\rho}P^{\ast}\right)  _{ij}  &  \leq\max
_{k}\left(  \mathrm{ord\,}p_{ik}+\mathrm{ord\,}\square^{t+\rho_{k}%
}+\mathrm{ord\,}p_{jk}\right) \\
&  \leq\max_{k}\left(  \sigma_{i}-\rho_{k}+\sigma_{j}-\rho_{k}+2\rho
_{k}+2t\right)  =\sigma_{i}+\sigma_{j}+2t
\end{align*}
The same inequality holds for the matrix $\square^{\sigma}Q^{\ast}%
\square^{t-\tau}Q\square^{\sigma}$ and for $\Omega.$ The principal symbol
$\mathrm{\Omega}\left(  z,\xi\right)  $ of $\Omega$ with respect to the shift
vector $2t+2\sigma$ is equal to%
\[
\mathrm{\Omega}=\mathrm{PR}^{t+\rho}\mathrm{P}^{\ast}+\mathrm{R}^{\sigma
}\mathrm{Q}^{\ast}\mathrm{R}^{t-\tau}\mathrm{QR}^{\sigma}%
\]
where $\mathrm{R}^{\omega}$ means the symbol of the operator $\square^{\omega
}.$ We will check that $\mathrm{\det\Omega}\left(  z,\xi\right)  \neq0$ as
$\xi\in\mathbb{R}^{n}\backslash\{0\}.$ If it is not the case for a point $\xi$
then there exists a non zero vector $v\in\mathbb{R}^{s}$ such that
$\mathrm{\Omega}\left(  z,\xi\right)  v=0.$ Define the coordinate scalar
product $\left\langle ,\right\rangle $ in $\mathbb{R}^{s}$ and write%
\begin{align*}
0  &  =\left\langle \mathrm{\Omega}\left(  x,\xi\right)  v,v\right\rangle
=\left\langle \mathrm{PR}^{t+\rho}\mathrm{P}^{\ast}v,v\right\rangle
+\left\langle \mathrm{R}^{\sigma}\mathrm{Q}^{\ast}\mathrm{R}^{t-\tau
}\mathrm{QR}^{\sigma}v,v\right\rangle \\
&  =\left\langle \mathrm{R}^{\left(  t+\rho\right)  /2}\mathrm{P}^{\ast
}v,\mathrm{R}^{\left(  t+\rho\right)  /2}\mathrm{P}^{\ast}v\right\rangle
+\left\langle \mathrm{R}^{\left(  t-\tau\right)  /2+\sigma}\mathrm{Q}%
v,\mathrm{R}^{\left(  t-\tau\right)  /2+\sigma}\mathrm{Q}v\right\rangle
\end{align*}
where $\mathrm{R}^{\omega/2}$ mean a diagonal matrix with the positive
diagonal terms\ $\sqrt{\mathrm{R}^{\omega_{i}}}$,$i=1,...,k.$ Both terms in
the right-hand side are non-negative, hence vanish. This yields
\begin{equation}
\mathrm{P}^{\ast}\left(  x,\xi\right)  v=0,\ \mathrm{Q}\left(  x,\xi\right)
v=0 \label{7}%
\end{equation}
By Proposition \ref{El} the sequence of symbols (\ref{6}) is exact at any real
point $\left(  x,\xi\right)  ,\xi\neq0$. Therefore the first equation
(\ref{7}) implies that $v=\mathrm{Q}^{\ast}\left(  x,\xi\right)  w$ for some
vector $w\in\mathbb{R}^{t}.$ By the second equation (\ref{7}) we find
$0=\left\langle \mathrm{Q}\left(  x,\xi\right)  v,w\right\rangle =\left\langle
v,v\right\rangle ,$ that is $v=0.$ This contradicts to the assumption and
completes the proof. $\blacktriangleright$

\begin{corollary}
\label{C}For any $x\in X$ and $r\leq1$ there exist linear continuous operators
$\mathrm{R}_{r}:\mathit{E}^{s}\left(  U_{x}\left(  r\right)  \right)
\rightarrow\mathit{E}^{r}\left(  U_{x}\left(  b_{x}^{3}r\right)  \right)  $
and $\mathrm{\Sigma}_{r}:\mathit{E}^{t}\left(  U_{x}\left(  r\right)  \right)
\rightarrow\mathit{E}^{s}\left(  U_{x}\left(  b_{x}^{3}r\right)  \right)  $
such that
\begin{equation}
\left(  P\mathrm{R}_{r}+\mathrm{\Sigma}_{r}Q\right)  g=g \label{10}%
\end{equation}
for any $g\in\mathit{E}^{s}\left(  U_{x}\left(  r\right)  \right)  $.
\end{corollary}

$\blacktriangleleft$ Write (\ref{4}) with two more terms%
\[
...\rightarrow\mathit{D}^{v}\overset{S}{\rightarrow}\mathit{D}^{u}\overset
{R}{\rightarrow}\mathit{D}^{t}\overset{Q}{\rightarrow}\mathit{D}^{s}%
\overset{P}{\rightarrow}\mathit{D}^{r}%
\]
and apply the functor $\mathrm{Hom}_{\mathit{D}}\left(  \cdot,\mathit{E}%
\right)  :$%
\[
\mathit{E}^{r}\overset{P}{\rightarrow}\mathit{E}^{s}\overset{Q}{\rightarrow
}\mathit{E}^{t}\overset{R}{\rightarrow}\mathit{E}^{u}\overset{S}{\rightarrow
}\mathit{E}^{v}\rightarrow...
\]
Here $P,Q,R,S,...$ are linear operators as in (\ref{5}). By Theorem \ref{cs}
applied to these terms, there exist linear continuous operators
\begin{align*}
\mathrm{\rho}_{r}  &  :\mathit{E}_{Q}\left(  U_{x}\left(  r\right)  \right)
\rightarrow\mathit{E}^{r}\left(  U_{x}\left(  b_{x}r\right)  \right)
,\mathrm{\sigma}_{r}:\mathit{E}_{R}\left(  U_{x}\left(  r\right)  \right)
\rightarrow\mathit{E}^{s}\left(  U_{x}\left(  b_{x}r\right)  \right)  \ \\
\mathrm{\tau}_{r}  &  :\mathit{E}_{S}\left(  U_{x}\left(  r\right)  \right)
\rightarrow\mathit{E}^{t}\left(  U_{x}\left(  b_{x}r\right)  \right)
\end{align*}
with the properties $P\mathrm{\rho}_{r}f=f$, $Q\mathrm{\sigma}_{r}g=g$,
$R\mathrm{\tau}_{r}h=h$ in $U_{x}\left(  b_{x}r\right)  $. We have $Q\left(
g-\mathrm{\sigma}_{r}Qg\right)  =0$ for any $g\in\mathit{E}^{s}\left(
U_{x}\left(  r\right)  \right)  .$ Therefore we can set $\mathrm{R}%
_{r}g=\mathrm{\rho}_{b^{2}r}\left(  g-\mathrm{\sigma}_{br}Qg\right)  $ and
similarly $\mathrm{\Sigma}_{r}h=\mathrm{\sigma}_{br}\left(  h-\mathrm{\tau
}_{r}Rh\right)  $. We have now for any $g\in\mathit{E}^{s}\left(  U_{x}\left(
r\right)  \right)  $%
\begin{align*}
\left(  P\mathrm{R}_{r}+\mathrm{\Sigma}_{r}Q\right)  g  &  =P\mathrm{\rho
}_{b^{2}r}\left(  g-\mathrm{\sigma}_{br}Qg\right)  +\mathrm{\sigma}%
_{br}\left(  Qg-\mathrm{\tau}_{r}RQg\right) \\
&  =g-\mathrm{\sigma}_{br}Qg+\mathrm{\sigma}_{br}Qg=g
\end{align*}
and (\ref{10}) follows. $\blacktriangleright$

\section{Solutions with compact support}

Let $U$ be an open set in $\mathbb{R}^{n};$ the topological dual space
$\mathit{E}^{\ast}\left(  U\right)  $ to $\mathit{E}\left(  U\right)  $ is
identified with space of distributions in $\mathbb{R}^{n}$ with compact
support contained in $U.$ An arbitrary differential $s\times r$-matrix $P$ in
$U$ with analytic coefficients defines a continuous operator $P:\mathit{E}%
^{r}\left(  U\right)  \rightarrow\mathit{E}^{s}\left(  U\right)  $ and the
adjoint operator $P^{\ast}:\mathit{E}^{\ast}\left(  U\right)  ^{s}%
\rightarrow\mathit{E}^{\ast}\left(  U\right)  ^{r}$ which acts by%
\[
P^{\ast}\phi=\psi,\ \psi\left(  u\right)  =\phi\left(  Pu\right)  ,\phi
\in\mathit{E}^{\ast}\left(  U\right)  ^{s},\ u\in\mathit{E}^{s}\left(
U\right)
\]
For any complex (\ref{4}) of left $\mathit{D}_{X}$-modules in an open set
$X\subset\mathbb{R}^{n}$ and any open set $U\subset X$ the sequence%
\[
...\rightarrow\mathit{E}^{\ast}\left(  U\right)  ^{t}\overset{Q^{\ast}%
}{\rightarrow}\mathit{E}^{\ast}\left(  U\right)  ^{s}\overset{P^{\ast}%
}{\rightarrow}\mathit{E}^{\ast}\left(  U\right)  ^{r}%
\]
is a complex of vector spaces.

\begin{theorem}
\label{ec}If (\ref{4}) is an elliptic complex in an open set $X\subset
\mathbb{R}^{n}$, then there exists a continuous function $c_{x}$ in $X$ such
that\newline\textbf{C}. for any point $x\in X$ and any $r,~0<r\leq1$ the
kernel of $P^{\ast}:\mathit{E}^{\ast}\left(  U_{x}\left(  c_{x}r\right)
\right)  ^{s}\rightarrow\mathit{E}^{\ast}\left(  U_{x}\left(  c_{x}r\right)
\right)  ^{r}$ is contained in the image of $Q^{\ast}:\mathit{E}^{\ast}\left(
U_{x}\left(  r\right)  \right)  ^{t}\rightarrow\mathit{E}^{\ast}\left(
U_{x}\left(  r\right)  \right)  ^{s}$,\newline\textbf{D}. a function
$\alpha\in\mathit{E}^{\ast}\left(  U_{x}\left(  c_{x}r\right)  \right)  ^{r}$
is equal to $P^{\ast}\beta$ for some $\beta\in\mathit{E}^{\ast}\left(
U_{x}\left(  r\right)  \right)  ^{s}$ if and only if $\alpha\left(  u\right)
=0$ for any $u\in\mathit{E}_{P}\left(  U_{x}\left(  c_{x}r\right)  \right)  .$
\end{theorem}

$\blacktriangleleft$ We omit the bottom index $x.$ Dualizing (\ref{10}) we get%

\[
\mathrm{R}_{r}^{\ast}P^{\ast}\alpha+Q^{\ast}\mathrm{\Sigma}_{r}^{\ast}%
\alpha=\alpha
\]
for an arbitrary $\alpha\in\mathit{E}^{\ast}\left(  U\left(  cr\right)
\right)  ^{s}.$ If $\alpha P=0,$ this equation yields $\alpha=Q\beta,$ where
$\beta=\mathrm{\Sigma}_{r}^{\ast}\alpha\in\mathit{E}^{\ast}\left(  U\left(
r\right)  \right)  ^{r}.$ This proves statement \textbf{C}.

Check \textbf{D}. If $\alpha=P^{\ast}\beta,$ then $u\left(  \alpha\right)
=Pu\left(  \beta\right)  =0.$ Vice versa, let $u\left(  \alpha\right)  =0$ for
any $u\in\mathit{E}_{P}\left(  U\left(  cr\right)  \right)  .$ The
distribution $\beta=\alpha-P^{\ast}\mathrm{R}_{r}^{\ast}\alpha$ fulfils
\thinspace$v\left(  \beta\right)  =u\left(  \alpha\right)  $ where
$w=v-\mathrm{R}_{r}Pv$ for an arbitrary $v\in\mathit{E}^{s}\left(  U\left(
r\right)  \right)  .$ We have $w\in\mathit{E}^{s}\left(  U\left(  cr\right)
\right)  $ and $Pw=\left(  P-P\mathrm{R}_{cr}P\right)  v=0$ since of
(\ref{10}). Therefore $w\left(  \alpha\right)  =0,$ hence $v\left(
\beta\right)  =0$ which yields $\beta=0$ and $\alpha=P^{\ast}\gamma,$
$\gamma=\mathrm{R}_{r}^{\ast}\alpha.$ $\blacktriangleright$

\section{Resolutions}

Now we take more invariant\ point of view on systems of equations like
(\ref{1}).

\textbf{Definition (\cite{M2}). }Let $\mathit{M}$ be a filtered
left\textit{\ }$\mathit{D}$\textit{-}module and%
\[
\mathrm{gr\,}\mathit{M}=\oplus_{k\in\mathbb{Z}}\mathit{M}_{k}/\mathit{M}_{k-1}%
\]
the corresponding graded $\mathrm{D}$-module. We assume that\newline%
(\textbf{i) }the \textrm{D}-module \textrm{gr\thinspace}$\mathit{M}$ is
finitely generated,\newline(\textbf{ii)} the $\mathit{O}$-module\textbf{\ }%
\textrm{gr\thinspace}$\mathit{M}_{x}$ is free.

\textbf{Definition. }Let $\mathit{M}$ be a left $\mathit{D}$-module satisfying
(\textbf{i).} The product \textrm{gr\thinspace}$\mathit{M}\otimes_{\mathit{O}%
}\mathbb{C}$ is a module of finite type over the polynomial algebra
$\mathrm{A}=\mathrm{D}\otimes\mathbb{C\cong C}\left[  \xi_{1},...,\xi
_{n}\right]  .$ The \textit{characteristic }set of $\mathit{M}$ is by
definition the support of $\mathit{V}=\mathit{V}\left(  \mathit{M}\right)
\mathit{\ }$in the support of the $\mathrm{A}$-module \textrm{gr\thinspace
}$\mathit{M}\otimes_{\mathit{O}}\mathbb{C}.$ The set $\mathit{V}$ is an
algebraic cone in the set $\mathbb{C}^{n}$ of maximal ideals of the algebra
$\mathrm{A}$. Any point $\xi\in\mathit{V}$ can be interpreted as a
multiplicative functional $\mu:\mathrm{gr\,}\mathit{M}\otimes_{\mathit{O}%
}\mathbb{C\rightarrow C}$ such that $\mu\left(  am\right)  =a\left(
\xi\right)  \mu\left(  m\right)  $ for arbitrary $a\in\mathrm{A,}%
m\in\mathrm{gr\,}\mathit{M}\otimes_{\mathit{O}}\mathbb{C}.$

We call $\mathit{M}$ \textit{elliptic} if the characteristic variety
$\mathit{V}\left(  \mathit{M}\right)  $ contains no real point $\xi\neq0.$

\textbf{Remark. }It is easy to check that the condition (\textbf{ii}) for
$\mathit{M}=\mathrm{Cok}P:\mathit{D}^{\sigma}\rightarrow\mathit{D}^{\rho}$ is
equivalent to (*) for $P$. The characteristic set $\mathit{V}$ of this module
coincide with set of points $\xi$ such that $\mathrm{rank}P\left(
x,\xi\right)  <s.$

\textbf{Definition.} Let$\ \alpha:\mathit{E}\rightarrow\mathit{F}$ be a
morphism of filtered \textit{D}-modules. It is called \textit{strict}, if it
agrees with the filtrations and $\alpha\left(  \mathit{E}_{k}\right)
=\alpha\left(  \mathit{E}\right)  \cap\mathit{F}_{k},\ k\in\mathbb{Z}.$

\begin{proposition}
\label{e}Let
\begin{equation}
\mathit{E}\overset{\alpha}{\rightarrow}\mathit{F}\overset{\beta}{\rightarrow
}\mathit{G} \label{24}%
\end{equation}
be a complex of morphisms of filtered vector spaces. If $\mathrm{Ker\,gr\,}%
\beta=\mathrm{\operatorname{Im}gr\,}\alpha$, the complex (\ref{24}) is exact
and $\alpha$ is strict.
\end{proposition}

$\blacktriangleleft$ Let $\beta\left(  f\right)  =0$ for an element
$f\in\mathit{F}.$ We have $f\in\mathit{F}_{k}$ for some $k$ and
\textrm{gr\thinspace}$\beta\left(  f\right)  =0.$ By the condition there is an
element $e_{0}\in\mathit{E}_{k}$ such that \textrm{gr\thinspace}$\alpha\left(
e_{0}\right)  =\mathrm{gr\,}f,$ that is $f-\alpha\left(  e_{0}\right)
\in\mathit{F}_{k-1}.$ The element $f^{\prime}\doteq f-\alpha\left(
e_{0}\right)  $ is contained in $\mathrm{Ker\,}\beta,$ we repeat the above
arguments with $k$ replaced by $k-1$ and obtain an element $e_{1}\in
\mathit{E}_{k-1}$ such that $f^{\prime}-\alpha\left(  e_{1}\right)
\in\mathit{F}_{k-2}$ and so on. Finally we get $f=\alpha\left(  e\right)  $
where $\alpha=\alpha_{0}+\alpha_{1}+...\in\mathit{E}_{k}.$
$\blacktriangleright$

\begin{proposition}
\label{gre}If the complex (\ref{24}) is exact, $\alpha$ and$\ \beta$ are
strict, then $\mathrm{Ker\,gr\,}\beta=\mathrm{\operatorname{Im}gr\,}\alpha.$
\end{proposition}

$\blacktriangleleft$ We have $\mathrm{gr\,}\beta\,\mathrm{gr\,}\alpha=0.$ Show
that $\mathrm{Ker\,gr\,}\beta\subset\mathrm{\operatorname{Im}gr\,}\alpha.$
Take an element $\mathrm{f}\in\mathrm{F}_{k}\doteq\mathit{F}_{k}%
/\mathit{F}_{k-1}$ such that \textrm{gr\thinspace}$\beta\left(  \mathrm{f}%
\right)  =0.$ Let $f\in\mathit{F}_{k}$ be an element of the class
$\mathrm{f.}$ We have $\beta\left(  f\right)  \in\mathit{G}_{k-1}$ and
$\beta\left(  f\right)  \in\beta\left(  \mathit{F}_{k-1}\right)  $ since
$\beta$ is strict, that is $\beta\left(  f-f^{\prime}\right)  =0$ for an
element $f^{\prime}\in\mathit{F}_{k-1}.$ We have $f-f^{\prime}=\alpha\left(
e\right)  ,e\in\mathit{E}$ since (\ref{24}) is exact and $\alpha\left(
e\right)  =\alpha\left(  e^{\prime}\right)  $ for some $e^{\prime}%
\in\mathit{E}_{k}$ since $\alpha$ is strict. This yields $\mathrm{f}%
=\mathrm{gr}\,\alpha\left(  \mathrm{e}\right)  $ where $\mathrm{e}$ is the
class of $e^{\prime}.$ $\blacktriangleright$

Let $\mathit{M}$ be a filtered left $\mathit{D}$-module and
\begin{equation}
...\rightarrow\mathit{D}^{\tau}\overset{Q}{\rightarrow}\mathit{D}^{\sigma
}\overset{P}{\rightarrow}\mathit{D}^{\rho}\overset{\pi}{\rightarrow}%
\mathit{M}\rightarrow0\label{21}%
\end{equation}
be a strict exact sequence of filtered left $\mathit{D}$-modules. The complex
of $\mathrm{D}$-modules%

\begin{equation}
...\rightarrow\mathrm{D}^{\tau}\overset{\mathrm{Q}}{\rightarrow}%
\mathrm{D}^{\sigma}\overset{\mathrm{P}}{\rightarrow}\mathrm{D}^{\rho}%
\overset{\mathrm{\pi}}{\rightarrow}\mathrm{gr}\,\mathit{M}\rightarrow0
\label{22}%
\end{equation}
is then well defined where all morphisms have degree $0.$ We call (\ref{22})
the \textit{principal part} of (\ref{21}).\ 

\begin{proposition}
\label{E}If a left $\mathit{D}$-module $\mathit{M}$ fulfils (\textbf{i)}, then
for\textbf{\ }any point $x\in X$ there exist a neighborhood $U$ of $x$ and a
resolution of the graded module $\mathrm{gr}\,\mathit{M}.$
\end{proposition}

$\blacktriangleleft\ $The product $\mathrm{gr\,}\mathit{M}\otimes_{\mathit{O}%
}\mathbb{C}$ is a module of the polynomial algebra $\mathrm{A}\doteq
\mathrm{D}\otimes\mathbb{C}.$ Construct a strict resolution of this module of
the form%
\begin{equation}
...\rightarrow\mathrm{D}^{\tau}\otimes\mathbb{C}\overset{\mathrm{Q}%
\otimes\mathbb{C}}{\rightarrow}\mathrm{D}^{\sigma}\otimes\mathbb{C}%
\overset{\mathrm{P}\otimes\mathbb{C}}{\rightarrow}\mathrm{D}^{\rho}%
\otimes\mathbb{C}\overset{\mathrm{\pi}}{\rightarrow}\mathrm{gr}\,\mathit{M}%
\otimes\mathbb{C}\rightarrow0 \label{32}%
\end{equation}
By (\textbf{i})$\ $there exists a surjective morphism $\mathrm{\pi}%
:\mathrm{D}^{r_{0}}\rightarrow\mathrm{gr\,}\mathit{M}.$ We choose a shift
vector $\rho=\left(  \rho_{1},...,\rho_{r_{0}}\right)  $ where $\rho_{i}%
\doteq\mathrm{ord\,\pi}_{\mathbb{C}}\left(  e_{i}\right)  ,~i=1,...,r_{0}$ for
the standard generators $e_{1},...,e_{r_{0}}$ and introduce the filtration
$\mathrm{D}^{\rho}$ in the module $\mathrm{D}^{r_{0}}$. The morphism
$\mathrm{\pi:D}^{\rho}\rightarrow\mathrm{gr\,}\mathit{M}$ has degree $0$ and
generate $\mathrm{A}$-morphism $\mathrm{\pi}_{\mathrm{A}}:\mathrm{A}^{\rho
}\rightarrow\mathrm{gr\,}\mathit{M}\otimes\mathbb{C}$ where $\mathrm{A}%
=\mathrm{D}\otimes\mathbb{C}.$ Because of the algebra $\mathrm{A}$ is
Noetherian,\textbf{\ }the submodule $\mathrm{Ker\,\pi}_{\mathrm{A}}$ is
generated by some homogeneous elements $\mathrm{p}_{1},...,\mathrm{p}_{r_{1}}%
$. Let $\mathrm{P}_{\mathrm{A}}:\mathrm{A}^{r_{1}}\rightarrow\mathrm{A}^{\rho
}$ be the morphism such that $\mathrm{P}_{\mathrm{A}}\left(  e_{i}^{\prime
}\right)  =\mathrm{p}_{i},$ $i=1,...,r_{1}$ for the standard generators
$e_{1}^{\prime},...,e_{r_{1}}^{\prime}$ of $\mathrm{A}^{r_{1}}.$ Set
$\sigma=\left(  \sigma_{1},...,\sigma_{r_{1}}\right)  $ where $\sigma
_{i}\doteq\mathrm{ord}_{\rho}\mathrm{p}_{i}$ and introduce the filtration
$\mathrm{D}^{\sigma}$ in\ $\mathrm{D}^{r_{1}}.$ The morphism $\mathrm{P}%
_{\mathrm{A}}$ is homogeneous of degree $0$ and $\mathrm{\operatorname{Im}%
\,P}_{\mathrm{A}}=\mathrm{Ker\,\pi}_{\mathrm{A}}\mathrm{.}$ Because of the
module $\mathrm{D}$ is Noetherian we can apply the same arguments to
$\mathrm{Ker\,P}_{\mathrm{A}}\ $and choose morphisms $\mathrm{Q}_{\mathrm{A}}$,....

By (\textbf{ii}) all the $n$ morphisms $\mathrm{P}_{\mathrm{A}},\mathrm{Q}%
_{\mathrm{A}},...$ have extensions to some $\mathrm{D}$-morphisms
$\mathrm{P},\mathrm{Q},...$ such that the sequence (\ref{22}) is a complex. It
can be shown by standard homological arguments since the $\mathit{O}$-modules
\textrm{gr}$\mathit{M},\mathrm{D}^{\rho},\mathrm{D}^{\sigma},...$ are flat.
$\blacktriangleright$

\begin{proposition}
\label{P}For any free graded resolution (\ref{22}) there exists a free strict
resolution (\ref{21}) of $\mathit{M}$ such that (\ref{22}) is the principal
part of \textit{(\ref{21}).}
\end{proposition}

$\blacktriangleleft$ Choose for any $i=1,...,r_{0}$ an element $m_{i}%
\in\mathit{M}_{\rho_{i}}$ whose image in $\mathit{M}_{\rho_{i}}/\mathit{M}%
_{\rho_{i}-1}$ is equal to $\mathrm{\pi}\left(  e_{i}\right)  $ and define a
$\mathit{D}$-morphism $\pi:\mathit{D}^{r_{0}}\rightarrow\mathit{M}$ such that
$\pi\left(  e_{i}\right)  =m_{i},i=1,...,r_{0}.$ This morphism agrees with the
filtrations in $\mathit{D}^{\rho}$ and in $\mathit{M}$ and is surjective,
since so is $\mathrm{\pi.}$ Next we lift $\mathrm{P}$ to a $\mathit{D}%
$-morphism $P_{0}:\mathit{D}^{r_{1}}\rightarrow\mathit{D}^{r_{0}}.$ For any
standard generator $e_{k}^{\prime}$ of \textrm{D}$^{r_{1}}$ the row
$\mathrm{p}_{k}\doteq\mathrm{P}\left(  e_{k}^{\prime}\right)  \in
\mathrm{D}^{\rho}$ satisfies $\mathrm{\pi p}_{k}=0,$ which means
$\pi\mathrm{p}_{k}\in\mathit{M}_{\rho_{k}-1},$ $k=1,...,r_{1}.$ Because of
exactness of (\ref{22}), there exists an element $q_{k}\in\mathit{D}^{r_{0}}$
such that \textrm{ord\thinspace}$q_{k}=\rho_{k}-1$ and $\pi\left(
q_{k}\right)  =\pi\mathrm{p}_{k}.$ We have $\pi\left(  \mathrm{p}_{k}%
-q_{k}\right)  \in\mathit{M}_{\rho_{k}-2}$ and so on up to filtration $-1$.
Finally we collect the lines $\mathrm{p}_{k}-q_{k}-q_{k}^{\prime}-...,$
$k=1,...,r_{1}$ in a matrix $P$ of size $r_{0}\times r_{1}$ and have $\pi
P=0.$ The principal part of the line $P\left(  e_{k}^{\prime}\right)  $ is
equal to $\mathrm{p}_{k},$ that is the principal part of\textrm{\ }$P$ is
$\mathrm{P}.$ By Proposition \ref{e} $P$ is strict and
$\mathrm{\operatorname{Im}\,}P=\mathrm{Ker\,}\pi.$

The image of the composition $P\mathrm{Q}:\mathit{D}^{r_{2}}\rightarrow
\mathit{D}^{r_{0}}$ is contained in $\mathrm{Ker\,}\pi$ and \textrm{ord}%
$_{\rho}P\mathrm{Q}\left(  e^{\prime\prime}\right)  <\mathrm{ord}_{\tau
}e^{\prime\prime}$ for each standard generator $e^{\prime\prime}$ of
$\mathit{D}^{r_{2}}.$ Because of $\mathrm{Ker\,\pi}=\mathrm{\operatorname{Im}%
\,P},$ there exists an element $q_{1}\in\mathit{D}^{r_{1}}$ such that
\textrm{ord}$_{\sigma}q_{1}=\mathrm{ord}_{\rho}P\mathrm{Q}\left(
e^{\prime\prime}\right)  \ $and $P\mathrm{Q}\left(  e^{\prime\prime}\right)
=\mathrm{P}q_{1}$ up to a term of filtration $<\mathrm{ord}_{\sigma}q_{1}%
.$\ We make a matrix $Q_{1}:\mathit{D}^{r_{2}}\rightarrow\mathit{D}^{r_{1}}$
from the lines $Q_{1}\left(  e^{\prime\prime}\right)  =q_{1}$ where
$e^{\prime\prime}$ runs over the set of generators of$\ \mathit{D}^{r_{2}}.$
Consider the composition $P\left(  \mathrm{Q}-Q_{1}\right)  :\mathit{D}%
^{r_{2}}\rightarrow\mathit{D}^{r_{0}}.$ We have now $\mathrm{ord}_{\rho
}P\left(  \mathrm{Q}-Q_{1}\right)  \left(  e^{\prime\prime}\right)
<\mathrm{ord}_{\rho}P\mathrm{Q}\left(  e^{\prime\prime}\right)  $ and can find
an element $q_{2}\in\mathit{D}^{r_{1}}$ such that \textrm{ord}$_{\sigma}%
q_{2}=\mathrm{ord}_{\rho}P\left(  \mathrm{Q}-Q_{1}\right)  \left(
e^{\prime\prime}\right)  $ up to a term of filtration $<\mathrm{ord\,}q_{2}.$
Define a matrix ~$Q_{2}$ by $Q_{2}\left(  e^{\prime\prime}\right)  =q_{2}$ for
the set of standard generators $e^{\prime\prime}$ then consider the matrix
$\mathrm{Q}-Q_{1}-Q_{2}$ and so on. This series is finite since
$...<\mathrm{ord}_{\sigma}q_{2}<\mathrm{ord}_{\sigma}q_{1}<\mathrm{ord}_{\tau
}e^{\prime\prime}.$ We set $Q=\mathrm{Q}-Q_{1}-Q_{2}-....$ By Proposition
\ref{e} $P_{1}$ is strict and $\mathrm{Im\,}Q=\mathrm{Ker\,}P.$ We construct a
matrix $R$ such that $\operatorname{Im}R=\mathrm{Ker}Q$ in the similar way and
so on.\ $\blacktriangleright$

\begin{proposition}
\label{El}If $\mathit{M}$ is an elliptic module then any strict resolution
resolution (\ref{21}) of $\mathit{M}$ is elliptic.
\end{proposition}

$\blacktriangleleft$ From (\ref{22}) we get the sequence%
\begin{equation}
...\rightarrow\mathrm{D}^{\tau}\otimes\mathbb{C}\overset{\mathrm{Q}%
\otimes\mathbb{C}}{\rightarrow}\mathrm{D}^{\sigma}\otimes\mathbb{C}%
\overset{\mathrm{P}\otimes\mathbb{C}}{\rightarrow}\mathrm{D}^{\rho}%
\otimes\mathbb{C}\rightarrow\mathrm{gr}\mathit{M}\otimes\mathbb{C}\rightarrow0
\label{31}%
\end{equation}
It is exact since of (\textbf{ii}). Take any real point $\xi\neq0;$ let
$\mathfrak{m}$ be the corresponding maximal ideal in $\mathrm{A.}$ Note that
all the terms$\ ...,\mathrm{D}^{\rho}\otimes\mathbb{C}$ are flat over
$\mathrm{A}$ and
\[
\mathrm{Tor}^{\ast}\left(  \mathrm{gr}\mathit{M}\otimes\mathbb{C}%
,\mathrm{A}\mathbf{/}\mathfrak{m}\right)  =0
\]
since $\left(  \mathrm{gr}\mathit{M}\otimes\mathbb{C}\right)  \otimes
_{\mathrm{A}}\mathrm{A}\mathbf{/}\mathfrak{m}=0$ because of $\xi$ does not
belong to the characteristic set of $\mathrm{gr}\mathit{M}\otimes
\mathbb{C}.\ $Therefore tensoring (\ref{31}) by $\mathrm{A}$-module
$\mathrm{A}/\mathfrak{m}$ we get the exact sequence%
\[
...\rightarrow\mathrm{D}^{\tau}\otimes\mathbb{C\otimes}_{\mathrm{A}}%
\mathrm{A}/\mathfrak{m}\overset{\mathrm{Q}\left(  x,\xi\right)  }{\rightarrow
}\mathrm{D}^{\sigma}\otimes\mathbb{C\otimes}_{\mathrm{A}}\mathrm{A}%
/\mathfrak{m}\overset{\mathrm{P}\left(  x,\xi\right)  }{\rightarrow}%
\mathrm{D}^{r}\otimes\mathbb{C\otimes}_{\mathrm{A}}\mathrm{A}/\mathfrak{m}%
\rightarrow0
\]
which proves the Proposition. $\blacktriangleright$

\section{Key Lemma}

Let $\left(  \mathit{M},\rho\right)  $ be a filtered $\mathit{D}$-module. The
set $\mathrm{Hom}_{\mathit{D}}\left(  \mathit{M},\mathit{D}\right)  $ of
$\mathit{D}$-morphisms $h:\mathit{M}\rightarrow\mathit{D}$ has a natural
structure of two-side $\mathit{D}$-module since $\mathit{D}$ has such a
structure. It possesses the dual filtration $\rho^{\ast}$ such that
\textrm{ord}$_{\rho^{\ast}}\left(  h\right)  =k$ if $h\left(  \mathit{M}%
_{i}\right)  \subset\mathit{D}_{i+k}$ for any $i.$ In particular
$\mathrm{Hom}\left(  \mathit{D}^{\rho},\mathit{D}\right)  \cong\mathit{D}%
^{-\rho},$ where $\mathit{D}^{-\rho}$ is a free $\mathit{D}$-module of the
same rank as $\mathit{D}^{\rho}$ with the shift vector $-\rho.$ Any morphism
of left $\mathit{D}$-modules $P:\mathit{D}^{\sigma}\rightarrow\mathit{D}%
^{\rho}$ generates the \textit{dual} morphism%
\[
P^{\prime}\doteq\mathrm{Hom}\left(  P,\mathit{D}\right)  :\mathit{D}^{-\rho
}\rightarrow\mathit{D}^{-\sigma},\ h\mapsto Ph
\]
where we interpret an element $h\in\mathit{D}^{-\rho}$ as a column. The map
$P^{\prime}$ is a morphism of \textit{right} $\mathit{D}$-modules.

Fix a point $x\in\mathbb{R}^{n}.$ Let
\begin{equation}
\mathit{R:}...\rightarrow\mathit{D}^{\rho_{2}}\overset{P_{1}}{\rightarrow
}\mathit{D}^{\rho_{1}}\overset{P_{0}}{\rightarrow}\mathit{D}^{\rho_{0}%
}\rightarrow0 \label{34}%
\end{equation}
be a strict resolution of a left $\mathit{D}$-module $\mathit{M}$ where
$\rho_{0},\rho_{1},\rho_{2},...\ $are some shift vectors. The complex
$\mathrm{Hom}_{\mathit{D}}\left(  \mathit{R},\mathit{D}\right)  $ looks as%
\begin{equation}
0\rightarrow\mathit{D}^{-\rho_{0}}\overset{P_{0}^{\prime}}{\rightarrow
}\mathit{D}^{-\rho_{1}}\overset{P_{1}^{\prime}}{\rightarrow}\mathit{D}%
^{-\rho_{2}}\rightarrow...\rightarrow\mathit{D}^{-\rho_{k-1}}\overset
{P_{k-1}^{\prime}}{\mathcal{\rightarrow}}\mathit{D}^{-\rho_{k}}\rightarrow...
\label{15}%
\end{equation}
where $\mathit{D}^{-\rho_{i}}$ is a free right $\mathit{D}$-module of the same
rank $r_{i}$ as $\mathit{D}^{\rho_{i}}$ and all morphisms agree with the
filtrations and $P^{\prime}$ means left multiplication of a column by a matrix
$P.$ It is a complex of right $\mathit{D}$-modules.

\begin{lemma}
\label{ex}If a left $\mathit{D}$-module $\mathit{M\ }$satisfies (\textbf{i,ii}%
),\ then the sequence (\ref{15}) is exact at the terms $\mathit{D}^{-\rho_{k}%
}$ with $k=0,...,m-1$ where $m=n-\mathrm{\dim}_{\mathbb{C}}\mathit{V}\left(
\mathit{M}\right)  .$
\end{lemma}

$\blacktriangleleft$ The principal part of (\ref{15}) is the complex of
modules over the graded commutative algebra $\mathrm{D}=\mathrm{D}%
_{x}\mathrm{:}$%
\[
\mathrm{\Pi}:0\rightarrow\mathrm{D}^{-\rho_{0}}\overset{\mathrm{P}_{0}%
^{\prime}}{\rightarrow}\mathrm{D}^{-\rho_{1}}\rightarrow...\rightarrow
\mathrm{D}^{-\rho_{k-1}}\overset{\mathrm{P}_{k-1}^{\prime}}%
{\mathcal{\rightarrow}}\mathrm{D}^{-\rho_{k}}\rightarrow...
\]
which is equal to $\mathrm{Hom}_{\mathrm{D}}\left(  \mathrm{R,D}\right)  $
where $\mathrm{R}$ is the principal part of (\ref{34}). We want to show that
this complex is exact in terms with $k=0,...,m-1.$ By Proposition \ref{gre}
$\mathrm{R}$ is a resolution of $\mathrm{gr\,}\mathit{M.}$ By the condition
(\textbf{ii) }the complex $\mathrm{R}\otimes\mathbb{C}$ is a resolution of
$\mathrm{gr\,}\mathit{M}\otimes\mathbb{C}$ over the algebra $\mathrm{A}%
\doteq\mathrm{D}\otimes\mathbb{C}$ where $\otimes=\otimes_{\mathit{O}_{x}}%
$which is isomorphic to the graded polynomial$\ \mathbb{C}$-algebra
$\mathbb{C}[\xi_{1},...,\xi_{n}].$ This yields
\begin{equation}
H^{k}\left(  \mathrm{Hom}_{\mathrm{A}}\left(  \mathrm{R}\otimes\mathbb{C}%
\mathrm{,A}\right)  \right)  \cong\mathrm{Ext}^{k}\left(  \mathrm{gr\,}%
\mathit{M}\otimes\mathbb{C},\mathrm{A}\right)  \label{17}%
\end{equation}
The right-hand side of (\ref{17}) vanishes for$\ k=0,...,m-1$ in virtue of
[\cite{P1},\ Corollary 1, \S 13] which means that the complex $\mathrm{Hom}%
_{\mathrm{A}}\left(  \mathrm{R}\otimes\mathbb{C}\mathrm{,A}\right)  ;$%
\begin{equation}
0\rightarrow\mathrm{D}_{x}^{-\rho_{0}}\otimes\mathbb{C}\overset{\mathrm{P}%
_{0}^{\prime}\otimes\mathbb{C}}{\mathbb{\rightarrow}}\mathrm{D}_{x}^{-\rho
_{1}}\otimes\mathbb{C\rightarrow}...\rightarrow\mathrm{D}_{x}^{-\rho_{k-1}%
}\otimes\mathbb{C}\overset{\mathrm{P}_{k-1}^{\prime}\otimes\mathbb{C}%
}{\rightarrow}\mathrm{D}_{x}^{-\rho_{k}}\otimes\mathbb{C}\rightarrow...
\label{23}%
\end{equation}
is acyclic in terms $\mathrm{D}_{x}^{-\rho_{k}}\otimes\mathbb{C},$
$k=0,...,m-1$. By definition for any shift vector $\omega$ we
have$\ \mathrm{D}^{\omega}\mathbb{=\oplus}_{i}\mathrm{D}_{i}^{\omega}$ where
$\mathrm{D}_{i}^{\omega}$ is the vector space of homogeneous elements of
grading $i$ and $\mathrm{D}^{\omega}\otimes\mathbb{C=\oplus}_{i}\mathrm{D}%
_{i}^{\omega}\otimes\mathbb{C}$ where $\mathrm{D}_{i}^{\omega}\otimes
\mathbb{C}$ is a finite dimensional space. The complex (\ref{23}) is
isomorphic to $\mathrm{\Pi}\otimes\mathbb{C}$ hence all morphisms are
homogeneous of degree zero. Therefore
\[
\mathrm{\Pi}=\oplus\mathrm{\Pi}_{i},\ \mathrm{\Pi}\otimes\mathbb{C=\oplus
}\mathrm{\Pi}_{i}\otimes\mathbb{C}%
\]
where $\mathrm{\Pi}_{i}$ is a complex of free $\mathrm{D}$-modules of finite
type and $\mathrm{\Pi}_{i}\otimes\mathbb{C\ }$is a finite dimensional
$\mathbb{C}$-vector spaces for each $i=0,1,....$ We know that the complex
$\mathrm{\Pi}_{i}\otimes\mathbb{C}$ is acyclic in degrees $<m.$ By Nakayama's
lemma the complex $\mathrm{\Pi}_{i}$ is acyclic in the same degrees for each
$i$ which implies that the same true for the complex $\mathrm{\Pi.}$ By
Proposition \ref{e} that (\ref{15}) is also acyclic in degrees $k=0,...,m-1$
and the morphisms $P_{0}^{\prime},...,P_{k-1}^{\prime}$ are strict.
$\blacktriangleright$

\begin{corollary}
\label{d}The complex (\ref{15}) is elliptic in any degree $k<m$.\ The
statement \textbf{C} of Theorem \ref{ec} holds for this complex and a
continuous function $d_{x}$ in $X.$
\end{corollary}

$\blacktriangleleft$ By Lemma \ref{ex} we need only to check that for any real
point $\xi\in\mathrm{SpectA}\cong\mathbb{C}^{n},\xi\neq0$ the complex%
\[
\mathrm{\Pi}_{x}\otimes\mathrm{A}/\mathfrak{m}%
\]
is acyclic where $\mathfrak{m\subset}\mathrm{A}$ is the maximal ideal of the
point\ $\xi.$ We have $H^{\ast}\left(  \mathrm{\Pi}_{x}\otimes\mathrm{A}%
/\mathfrak{m}\right)  \cong H^{\ast}\left(  \mathrm{\Pi}_{x}\right)
\otimes\mathrm{A}/\mathfrak{m}$ since $\mathrm{A}$-module $\mathrm{\Pi}_{x}$
is flat. By Lemma \ref{ex} $H^{\ast}\left(  \mathrm{\Pi}_{x}\right)  =0.\ $The
second statement follows from Theorem \ref{ec}. $\blacktriangleright$

\section{Extension of solutions of determined systems}

\textbf{Definition. }Let $\mathit{M}$ be a left $\mathit{D}$-module with good
filtration that fulfils the condition (\textbf{i}). The characteristic set
$\mathit{V}=\mathit{V}\left(  \mathrm{gr}\mathit{M}\right)  $ is an algebraic
cone in $\mathbb{C}^{n}.$ We say that$\ \mathit{M}\mathcal{\ }$is
\textit{underdetermined} in $\mathit{V}=\mathbb{C},$ \textit{determined}, if
$\mathrm{\dim}_{\mathbb{C}}\mathit{V}<n$ and \textit{overdetermined}, if
$\mathrm{\dim}_{\mathbb{C}}\mathit{V}<n-1.$

Let $U$ be an open set in $X$ and $K\Subset U.$ If $\mathit{M}$ is not a
overdetermined module, then a solution $u$ of $\mathit{M}$ in a domain
$U\backslash K$ may have nonremovable singularity on $K,$ e.g. a \ fundamental
solution of the operator $P.$ A necessary condition for a solution $u$ to have
an extension to $U$ as a solution is vanishing of some momenta. Fix a smooth
density $\phi$ with support in an open set $V\subset U\ $such that
$\phi=\mathrm{d}x$ in a neighborhood of $K,$ take an arbitrary solution $v$ of
$P^{\ast}v=0$ defined in a neighborhood $V$ of $K$ such that $\mathrm{supp\,}%
\nabla\phi\Subset V\backslash K\ $and consider the integral%
\begin{equation}
\int_{U\backslash K}uP^{\ast}\left(  \phi v\right)  \label{11}%
\end{equation}
Note that if $u$ has extension to $U$ as solution of $\mathit{M,}$ we can
integrate in (\ref{11}) by parts and get the equation $\int\phi vP\left(
u\right)  =0.$ We state an inverse implication:

\begin{theorem}
\label{K}Let $\mathit{M}$ be an elliptic $\mathit{D}$-module in $X$ and
$c_{x}$ be the function in $X$ as in Theorem \ref{ec}. Let $x\in X,$ $0<r\leq
c_{x},$ $U_{x}=U_{x}\left(  r\right)  ,V_{x}=U_{x}\left(  c_{x}r\right)  $ and
$K\subset V_{x}$ be a compact set without holes. Then an arbitrary solution of
$Pu=0$ defined in $U_{x}\backslash K$ has a unique extension to $U_{x}$ as a
solution provided the integral (\ref{11}) vanishes for any smooth solution $v$
of $P^{\ast}v=0$ in $V_{x}.$
\end{theorem}

$\blacktriangleleft$ We may assume that $\mathrm{supp\,}\phi\subset V_{x}$ and
$\phi\left(  x\right)  =\mathrm{d}x$ \ in a neighborhood $W$ of $K.$ Take a
smooth function $e$ in $\mathbb{R}^{n}$ supported in $W$ that is equal to $1$
in a neighborhood $W_{0}$ of $K.$ The function $P\left(  eu\right)  $ is
supported in $W$ and vanishes in $W_{0}.$ Set $\alpha=P\left(  eu\right)
\mathrm{d}x$ in $V_{x}\backslash K$ and $\alpha=0$ in $K.$ We have $\alpha
\in\mathit{E}^{\ast}\left(  V_{x}\right)  ^{s}$ and for any solution $w$ of
the equation $P^{\prime}w=0$ in $V_{x}$
\[
\alpha\left(  w\right)  =\int_{V_{x}\backslash K}P\left(  eu\right)
w\mathrm{d}x=\int euP^{\ast}\left(  w\mathrm{d}x\right)  =\int uP^{\ast
}\left(  \phi w\mathrm{d}x\right)
\]
since the distribution $P^{\ast}\left(  \phi w\mathrm{d}x\right)  $ is
supported in $V_{x}\backslash W.$ By the assumption the right-hand side
vanishes for any $w.$ By Theorem \ref{ec} \textbf{D }there exists a
distribution $\beta\in\mathit{E}^{\ast}\left(  U_{x}\right)  ^{r}$ such that
$P\beta=\alpha,$ that is $Pu^{\prime}=0$ in $U_{x}$ where $u^{\prime}\doteq
eu-\beta.$ The functions $u$ and $u^{\prime}$ coincide in $U_{x}%
\backslash\mathrm{supp\,}e\cup\mathrm{supp\,}\beta$ and are analytic in
$U_{x}\backslash K,$ hence $u^{\prime}=u$ in $K$ since $K$ has no holes.
$\blacktriangleright$

\section{Overdetermined systems}

Let $P\ $be a matrix differential operator as in (\ref{1}) with analytic
coefficients defined in an open set $X\subset\mathbb{R}^{n}.$ We assign to
this matrix a sheaf of differential modules in $X.$ For this we globalize the
construction of Section 7: let$\ \mathcal{O}$ be the sheaf of germs of
analytic functions in $X$ and $\mathcal{D}$ be the sheaf algebra of germs of
differential operators with coefficients in $\mathcal{O}$. The stalk of
$\mathcal{D}$ at a point $x\in X$ is the algebra $\mathit{D}_{x}$ as in
Section 2. Let $\mathcal{M}$ be a filtered left $\mathcal{D}$-module defied in
open set $X\subset\mathbb{R}^{n}\ $that can be included in a strict exact
sequence of filtered left $\mathcal{D}$-modules%
\begin{equation}
\mathcal{D}^{\sigma}\overset{P}{\rightarrow}\mathcal{D}^{\rho}\overset{\pi
}{\rightarrow}\mathcal{M}\rightarrow0\label{37}%
\end{equation}
where $\mathcal{D}^{\sigma},\mathcal{D}^{\rho}$ denote some filtrations in
free left $\mathcal{D}$-modules defined as in Section 3.\ Here $P$ acts as a
morphism of left $\mathcal{D}$-modules: $a\mapsto aP\ $and the filtration in
$\mathcal{M}$ is the image of the filtration in $\mathcal{D}^{\rho}:$
$\mathcal{M}_{k}=\pi\left(  \mathcal{D}_{k}^{\rho}\right)  ,k\in Z.$ Note that
for any function $\mathcal{O}$-sheaf $\mathcal{F\ }$in $X$ the space
$\mathrm{Hom}_{\mathcal{O}}\left(  \mathcal{M},\mathcal{F}\right)  $ is
isomorphic to the space of solutions of the equation (\ref{1}) in the space
$\Gamma\left(  X,\mathcal{F}\right)  .$

\begin{proposition}
For any compact set $K\subset X,$ the sequence (\ref{37}) can be completed to
a strict exact complex of $\mathcal{D}$-sheaves%
\begin{equation}
...\rightarrow\mathcal{D}^{\tau}\overset{Q}{\rightarrow}\mathcal{D}^{\sigma
}\overset{P}{\rightarrow}\mathcal{D}^{\rho}\overset{\pi}{\rightarrow
}\mathcal{M}\rightarrow0\label{38}%
\end{equation}
defined in a neighborhood of $K$ where $...,\mathcal{D}^{\tau}$ are filtered
free $\mathcal{D}$-sheaves of the same type.
\end{proposition}

$\blacktriangleleft$ Let \textrm{D}$_{X}$ be the sheaf in $X$ whose stalk are
the algebras \textrm{D }and\textrm{ }$\mathrm{D}_{X}^{\omega}\ $be the graded
\textrm{D}$_{X}$-sheaf where $\omega$ is an arbitrary shift vector. Consider
the sequence of graded $\mathrm{D}_{X}$-modules%
\[
\mathrm{D}_{X}^{\sigma}\overset{\mathrm{P}}{\rightarrow}\mathrm{D}_{X}^{\rho
}\overset{\mathrm{\pi}}{\rightarrow}\mathrm{gr}\mathcal{M}\rightarrow0
\]
generated by (\ref{37}). It is exact since of Proposition\ \ref{gre}. For
$k=0,1,2,...$ we consider $\mathcal{O}$-sheaf $\left(  \mathrm{Ker~P}\right)
_{k}:\left(  \mathrm{D}_{X}^{\sigma}\right)  _{k}\rightarrow\left(
\mathrm{D}_{X}^{\rho}\right)  _{k}.$ It is a coherent analytic sheaf in the
real domain $X.$ Let $L$ be a compact set in $X$ such that $K\Subset L.$ By
the classical theory of coherent sheaves the sheaf $\left(  \mathrm{Ker~P}%
\right)  _{k}$ is generated in each point $x\in L$ by a finite set
\textrm{S}$_{k}$ of its sections. The total set $\mathrm{S}\doteq\cup
_{k}\mathrm{S}_{k}$ generates $\mathrm{D}_{x}$-sheaf \textrm{Ker\ P}%
$_{x}:\mathrm{D}_{x}^{\sigma}\rightarrow\mathrm{D}_{x}^{\rho}$ at each point
$x\in L.$ On the other hand, for any point $x$ there is a finite subset
$\mathrm{q}_{x}\subset\mathrm{S}$ that generates the stalk $\left(
\mathrm{Ker\ P}\right)  _{x}$ since the algebra \textrm{D}$_{x}$ is
Noetherian. Obviously the set $\mathrm{q}_{x}$ generates the sheaf
\textrm{Ker~P} also in a neighborhood of $x.$ Therefore there is a finite set
$F\subset L$ such that the union \textrm{q}$_{L}=\cup\left\{  \mathrm{q}%
_{x},x\in F\right\}  $ generates the $\mathrm{D}$-sheaf \textrm{Ker P} at each
point $x\in L.$ Let $\mathrm{D}_{X}^{t}$ be a free \textrm{D}-sheaf with
generators $e_{1},...,e_{t}.$ Consider a \textrm{D}$_{X}$-morphism
\textrm{Q}$:\mathrm{D}_{X}^{t}\rightarrow\mathrm{D}_{X}^{\sigma}$ such that
$\mathrm{q}_{j}=\mathrm{Q}\left(  e_{j}\right)  ,\ j=1,...,t$\ are all
elements of the set \textrm{q}$_{L}$. Define a filtration $\mathrm{D}%
_{X}^{\tau}$ in \textrm{D}$_{X}^{t}$ by means of a shift vector $\tau=\left(
\tau_{1},...,\tau_{t}\right)  $ where $\tau_{j}=\mathrm{\deg q}_{j}%
,j=1,...,t.$ The morphism $\mathrm{Q}:$ \textrm{D}$^{\tau}\rightarrow
\mathrm{D}^{\rho}$\ agrees with the filtrations and $\operatorname{Im}%
\mathrm{Q}=\mathrm{Ker~P}$.\ Next we consider the restriction \textrm{Q}$_{L}$
of \textrm{Q}$_{X}$ to $L$ and repeat these arguments for the $\mathrm{D}%
$-sheaf $\mathrm{Ker~Q}_{L}$ and so on. We obtain in this way an exact
sequence of \textrm{D}$_{Y}$-sheaves%

\[
...\overset{\mathrm{R}_{x}}{\rightarrow}\mathrm{D}_{Y}^{\tau}\overset
{\mathrm{Q}_{x}}{\rightarrow}\mathrm{D}_{Y}^{\sigma}\overset{\mathrm{P}%
}{\rightarrow}\mathrm{D}_{Y}^{\rho}\overset{\mathrm{\pi}}{\rightarrow
}\mathrm{gr}\mathcal{M}_{Y}\rightarrow0
\]
defined in a neighborhood $Y$ of $K.$ Then we construct a strict exact
sequence (\ref{38}) by means of arguments of Propositions \ref{E}%
.\ $\blacktriangleright$

Note that for any point $x\in X$ the stalks of (\ref{38}) form a resolution of
$\mathit{D}$-module $\mathcal{M}_{x}$ like (\ref{21}).

\begin{theorem}
\label{ee}Suppose that a filtered left $\mathcal{D}$-module $\mathcal{M\ }$in
$X$ can be included in a strict exact sequence of the form (\ref{37}) such
that the \textrm{D}$_{x}$-module \textrm{gr}$\mathcal{M}_{x}$ fulfils the
condition (\textbf{ii}), is elliptic and overdetermined at any point $x\in X.$
Let $Y$ be a relatively compact subset of $X$ and $S$ be a closed $C^{1}%
$-submanifold of $Y$ of dimension $d=n-\max_{X}\mathrm{\dim}_{\mathbb{C}%
}\mathit{V}\left(  \mathrm{gr}\,\mathcal{M}_{x}\right)  -2.$ There exists an
open neighborhood $V$ of $S$ such that any solution $u$ of \textit{(\ref{1})}
in $Y\backslash\bar{V}$ has a unique extension in $X$ as a solution of
\textit{the same system.}
\end{theorem}

\textbf{Remark.} If $X$ is connected the function $x\mapsto\mathrm{\dim
}_{\mathbb{C}}\mathit{V}\left(  \mathrm{gr}\,\mathcal{M}_{x}\otimes
\mathbb{C}\right)  $ is constant in $X$. This follows from (\textbf{ii}).

\textbf{Example 1.} Let $d=0,\ $then the statement tells that for any point
$x\in X$ there exists a compact neighborhood $K_{x}\subset X$ such that
arbitrary solution defined in $X\setminus K$ is uniquely extended to a
solution in $X.$

$\blacktriangleleft$ Proof of Theorem. Introduce an Euclidean structure in
$\mathbb{R}^{n}$.

\begin{lemma}
\label{le}There exist positive constants $b\leq c<1$ that depends only on $K$
such that for\ an arbitrary subspace $Z$ in $\mathbb{R}^{n}$ of dimension
$d\ $and arbitrary balls $Y\left(  r\right)  $ and $Z\left(  s\right)  $ of
radius $r$ in $Y=Z^{\perp},$ respectively in $Z$ of radii $r,s\leq1,\ cr\leq
s$ such that $Y\left(  r\right)  \times Z\left(  s\right)  \subset K$ an
arbitrary solution $u$ of $\mathcal{M}$ defined in the set $Y\left(  r\right)
\backslash Y\left(  br\right)  \times Z\left(  s\right)  $ has a unique
extension to $Y\left(  r\right)  \times Z\left(  s\right)  $ as a solution of
$\mathit{M}.$
\end{lemma}

Here and later we denote by $Y\left(  r^{\prime}\right)  $ the ball with the
same center as $Y\left(  r\right)  $; notation $Z\left(  s^{\prime}\right)
\ $has a similar meaning. To prove the Theorem we take for an arbitrary point
$x_{0}\in S$ the tangent subspace $Z$ to $S$ at $x_{0}$ and set $Y=Z^{\perp}.$
In the case $d=0$ we take $Z=0,Y=\mathbb{R}^{n}.$ Choose a positive number $r$
such that $Y\left(  r\right)  \backslash Y\left(  br\right)  \times Z\left(
cr\right)  \subset X\backslash S.$ This choice is possible since $S$ is
contained in $o\left(  r\right)  $-neighborhood of $Z$. By Lemma \ref{le} any
solution $u$ can be extended to the set $Y\left(  r\right)  \times Z\left(
cr\right)  .$ This set contains a neighborhood of $x_{0}.$ We take for $V$ the
union of these neighborhoods for all $x_{0}\in S$ and complete the proof of
Theorem. $\blacktriangleright$

$\blacktriangleleft$ Proof of Lemma\ \ref{le}. Choose some positive numbers
$r_{0},s_{0}\leq1$ such that\ $Y\left(  r_{0}\right)  \times Z\left(
s_{0}\right)  \Subset X;$ we may assume that $r_{0}=1,s_{0}=c$ by coordinate
rescaling. Set $b=c^{m+1},\ c=\inf_{U}c_{x}/4\ $where $c_{x}$ is the function
as in Lemma \ref{d}. Choose a coordinate system $\left(  y,z\right)  $ in
$Y\times Z$ such that the centers of $Y\left(  r\right)  $ and $Z\left(
r\right)  $ are in the origins.

Take a smooth function $e$ in $Y$ with support in $Y\left(  2b\right)  $ such
that $e=1$ in $Y\left(  b+\varepsilon\right)  $ and set $v\left(  x\right)
=P_{0}\left(  e\left(  y\right)  u\left(  x\right)  \right)  ,$ the function
$v_{0}$ is extended by zero to $Y\left(  b\right)  \times Z\left(  c\right)
.$ Take a convex polytope $\Pi\subset Z\left(  c\right)  \backslash Z\left(
c/2\right)  ;$ let $F_{\alpha},\alpha\in N$ be its faces. Let $N_{k}$ be the
subset of $N$ of faces $F_{\alpha}$ of dimension $k=0,1,...,\dim Z;\ $the
face$\ \Pi$ is the only one of dimension $d\doteq\dim Z.$ The notation
$\alpha_{k}$ always will mean that $\alpha_{k}\in N_{k}$. We suppose that each
face $F_{\alpha_{k}}$ of dimension $k<d$ is a simplex and the inequality
holds
\begin{equation}
2b\leq\mathrm{diam\,}F_{\alpha_{1}}\leq3b \label{28}%
\end{equation}
for each 1-face. We call $k$-\textit{flag} any sequence $A=\left(  \alpha
_{k},\alpha_{k+1},...,\alpha_{d-1}\right)  $ such that $F_{\alpha_{k}}\subset
F_{\alpha_{k+1}}\subset...$.$\subset F_{\alpha_{d-1}}.$ For a set $G\subset Z$
and a positive $\varepsilon$ we denote by $\left(  G\right)  _{\varepsilon}$
the open $\varepsilon$-neighborhood of $G.$

Take a smooth function $f_{0}$ in $Z$ with compact support in $\left(
\Pi\right)  _{b}\ $such that $f_{0}=1$ in $\Pi.$ For an arbitrary $k<d$ and
$\alpha_{k}\in N_{k}$ we choose a smooth function $f_{\alpha_{k}}$ that
fulfils\newline\textbf{I.} \textrm{supp}$\mathrm{\,}f_{\alpha_{k}}%
\subset\left(  F_{\alpha_{k}}\right)  _{b}$ and\newline\textbf{II.}
$\sum_{\alpha_{k}}f_{\alpha_{k}}=1$ in $\left(  \cup_{\alpha_{k}}F_{\alpha
_{k}}\right)  _{b}.$\newline Take an arbitrary $k$-flag $\mathit{A}=\left(
\alpha_{k},\alpha_{k+1},...\right)  $ and define the function%
\begin{equation}
v_{\mathit{A}}=P_{d-k+1}\left(  f_{\alpha_{k}}...P_{2}\left(  f_{\alpha_{d-1}%
}P_{1}\left(  f_{0}v_{0}\right)  \right)  ...\right)  \label{12}%
\end{equation}
where $P_{1},...,P_{d+1}$ are differential operators as in (\ref{15}) (strokes
are omitted).

\begin{lemma}
\textbf{III.} The function $v_{\mathit{A}}$ is supported by $\left(
F_{\alpha_{k}}\right)  _{b}$.\newline\textbf{IV.} For any $k+1$-flag
$\mathit{B}$ we have%
\[
\sum_{\alpha_{k}}v_{\alpha_{k},\mathit{B}}=0
\]
where the sum is taken for all $k$-flags that contain $\mathit{B}.$\newline
\end{lemma}

$\blacktriangleleft$ Statement \textbf{III} follows from \textbf{I} and
equation \textbf{IV} follows from \textbf{II}:%
\[
\sum_{\alpha_{k}}v_{\alpha_{k},\mathit{B}}=P_{d-k+1}\sum_{\alpha_{k}}%
f_{\alpha_{k}}v_{\mathit{B}}=P_{d-k+1}v_{\mathit{B}}=0\ \blacktriangleright
\]

For any $1$-flag $\mathit{A}=\left(  \alpha_{1},...\right)  $ we have
$v_{\mathit{A}}=v_{\alpha_{0},\mathit{A}}+v_{\beta_{0},\mathit{A}}$ where
$\alpha_{0},\beta_{0}\in N_{0}$ are the vertices of the face $F_{\alpha_{1}}$
hence $\left(  \alpha_{0},\mathit{A}\right)  $ and $\left(  \beta
_{0},\mathit{A}\right)  $ are $0$-flags. By \textbf{III} we have
\textrm{supp}$\mathrm{\,}v_{\alpha_{0},\mathit{A}}\Subset\left(  F_{\alpha
_{0}}\right)  _{b}$ and similarly for $v_{\beta_{0},\mathit{A}}.$ the left
inequality (\ref{28}) implies that the supports of the distributions
$v_{\alpha_{0},\mathit{A}}$ and $v_{\beta_{0},\mathit{A}}$ are disjoint. The
formula (\ref{12}) yields $P_{d+1}v_{\alpha_{0},\mathit{A}}=P_{d+1}%
v_{\beta_{0},\mathit{A}}=0$ hence by Lemma \ref{d} there exist solutions to
the equations%
\begin{equation}
v_{\alpha_{0},\mathit{A}}=P_{d}w_{\alpha_{0},\mathit{A}},\ v_{\beta
_{0},\mathit{A}}=P_{d}w_{\beta_{0},\mathit{A}} \label{30}%
\end{equation}
with compact supports \textrm{supp}$\mathrm{\,}w_{\alpha_{0},\mathit{A}%
}\Subset\left(  F_{\alpha_{0}}\right)  _{b/2c},$ \textrm{supp}$\mathrm{\,}%
w_{\beta_{0},\mathit{A}}\Subset\left(  F_{\beta_{0}}\right)  _{b/2c}.$ Set
$w_{\mathit{A}}=w_{\alpha_{0},\mathit{A}}+w_{\beta_{0},\mathit{A}}$ and have
$P_{d}w_{\mathit{A}}=v_{\mathit{A}}.$ By (\ref{28}) for any $\alpha_{0}%
,\ $\textrm{supp}$\mathrm{\,}w_{\alpha_{0},\mathit{A}}\Subset\left(
F_{\alpha_{1}}\right)  _{b/2c}\subset\left(  F\right)  _{3b+b/2c}%
\subset\left(  F\right)  _{b/c}$ since $3b+b/2c\leq b/c.$ By \textbf{IV} we
have
\[
\sum_{\alpha_{0},\alpha_{1}}v_{\alpha_{0},\alpha_{1},\mathit{B}}=\sum
_{\alpha_{1}}v_{\alpha_{1},\mathit{B}}=0
\]
where the sum is taken over all flags that contain the $2$-flag $\mathit{B}%
=\left(  \alpha_{2},\alpha_{3},...\right)  .$ Therefore we can assume that
also
\begin{equation}
\sum_{\alpha_{0},\alpha_{1}}w_{\alpha_{0},\alpha_{1},\mathit{B}}=\sum
_{\alpha_{1}}w_{\alpha_{1},\mathit{B}}=0 \label{29}%
\end{equation}
Define $v_{\mathit{A}}^{\prime}=v_{\mathit{A}}-\sum w_{\beta,\mathit{A}}$ for
any $1$-flag $\mathit{A}.$ By (\ref{28}) we have $\mathrm{supp\,}%
v_{\mathit{A}}^{\prime}\Subset\left(  F_{\alpha_{0}}\right)  _{b/c}\ $for any
$1$-flag $\mathit{A}$ and an arbitrary vertex $F_{\alpha_{0}}$ of the face
$F_{\alpha_{1}}.$ Due to (\ref{30}) we have $P_{d}v_{\mathit{A}}^{\prime}=0,$
hence by Lemma \ref{d} there exists a solution $w_{\mathit{A}}$ to
$P_{d-1}w_{\mathit{A}}=v_{\mathit{A}}^{\prime}$ with compact support in
$\left(  F_{\alpha_{0}}\right)  _{b/2c^{2}}.$ Set for any $2$-flag
$\mathit{B}$%
\[
v_{\mathit{B}}^{\prime}=v_{\mathit{B}}-\sum_{\alpha_{1}}w_{\alpha
_{1},\mathit{B}}%
\]
where the sum is taken over all $1$-flags that contains the flag $\mathit{B}$.
We have $\mathrm{supp\,}v_{\mathit{B}}^{\prime}\Subset\left(  F_{\alpha_{0}%
}\right)  _{b/c^{2}}.$ By (\ref{29}) and (4) we have%
\[
P_{d-1}v_{\mathit{B}}^{\prime}=P_{d-1}\sum_{\alpha_{1}}\left(  f_{\alpha_{1}%
}v_{\mathit{B}}-w_{\alpha_{1},\mathit{B}}\right)  =\sum_{\mathit{B}%
\subset\mathit{A}}\left(  v_{\mathit{A}}-P_{d-1}w_{\mathit{A}}\right)
=\sum_{\alpha_{1}}w_{\alpha_{1},\mathit{B}}=0
\]
By Lemma \ref{d} we can solve the equation $P_{d-2}w_{\mathit{B}%
}=v_{\mathit{B}}^{\prime}$ for a function $w_{\mathit{B}}^{\prime}$ with
compact support in $\left(  F_{\alpha_{0}}\right)  _{b/2c^{3}}$ under the
condition
\[
\sum_{\alpha_{2}}w_{a_{2},\mathit{C}}=0
\]
for any $3$-flag $\mathit{C}.$ Set
\[
v_{\mathit{C}}^{\prime}=v_{\mathit{C}}-\sum w_{\alpha_{2},\mathit{C}}%
\]
and have $P_{d-2}v_{\mathit{C}}^{\prime}=0$ for any $3$-flag $\mathit{C.}$
Continuing arguing in this way $d-1$ time, we get the function%
\[
v_{0}^{\prime}=v_{0}-\sum_{\alpha_{d-1}}w_{\alpha_{d-1}}%
\]
where $\mathrm{supp\,}w_{\alpha_{d-1}}\Subset\left(  F_{\alpha_{0}}\right)
_{b/c^{d-1}}$ and $P_{1}v_{0}^{\prime}=0.$ We have $\mathrm{supp\,}%
v_{0}^{\prime}\Subset Y\left(  2c\right)  \times Z\left(  c\right)  $ and
$v_{0}^{\prime}=v_{0}$ in $Y\times Z\left(  c\right)  .$ We apply again Lemma
\ref{d} and find a solution to the equation $P_{0}w_{0}=v_{0}^{\prime}$ with
compact support in $Y\left(  1/2\right)  \times Z\left(  1/2\right)  .$ We
have $P_{0}w_{0}=f_{0}P_{0}\left(  eu\right)  $ in $Y\times Z\left(  c\right)
$ Because of $f_{0}=1$ in $Y\times\Pi$ we have $P_{0}\left(  eu-w_{0}\right)
=0,$ hence $P_{0}\left(  \left(  1-e\right)  u+w_{0}\right)  =0$ in
$Y\times\Pi.$ The function $U\doteq\left(  1-e\right)  u+w_{0}$ fulfils the
equation $P_{0}U=0$ and coincides with $u$ in $Y\left(  1\right)  \backslash
Y\left(  1/2\right)  \times Z\left(  c-\varepsilon\right)  .$ By uniqueness of
analytic continuation we have $U=u$ in $Y\left(  1\right)  \backslash Y\left(
c\right)  \times Z\left(  c-\varepsilon\right)  $ for arbitrary $\varepsilon
>0.\blacktriangleright$

\textbf{Example 2. }The statement of Theorem \ref{ee} does not hold in general
for $\dim\,S=\mathrm{\dim}_{X}\mathit{V}\left(  \mathrm{gr}\,\mathit{M}%
\right)  +1.\ $Let $\mathbb{R}^{n}=Y\oplus Z,$ where $Z$ is spanned by the
coordinates $x_{1},...,x_{d}.$ Consider the $\mathit{D}$-module $\mathit{M}%
=\mathit{D}/\left(  p_{0},p_{1},...,p_{d}\right)  $ where%
\[
p_{0}=p_{0}\left(  \partial_{x_{d+1}},...,\partial_{x_{n}}\right)
,\ p_{i}=\partial_{x_{i}},i=1,...,d
\]
where $p_{0}$ is an elliptic operator with constant coefficients in $Y.$ It is
an elliptic module and $\mathit{\ V}\left(  \mathrm{gr}\,\mathcal{M}%
_{x}\right)  =\{\left(  x,\xi\right)  ;\xi_{1}=...=\xi_{d}=\mathrm{p}%
_{0}\left(  \xi_{d+1},...,\xi_{n}\right)  =0\}$ for any $x\in\mathbb{R}^{n}.$
Dimension of this characteristic manifold is equal to$\ n-d-1$ however there
is no compulsory extension for solutions of $\mathit{M}$ from $\mathbb{R}%
^{n}\backslash Z$ on $\mathbb{R}^{n}$ since any fundamental solution $E$ of
$p_{0}$ considered as a function in $\mathbb{R}^{n}$ has singularity in $Z.$

\end{document}